\documentclass[a4paper, 10pt]{article}
\usepackage{graphics}
\usepackage{color}
\usepackage[pdfencoding=auto,psdextra]{hyperref}

\usepackage{lscape}
\usepackage{latexsym}
\usepackage{amsmath,verbatim}
\usepackage{amsthm}
\usepackage{amssymb}
\usepackage[mathscr]{euscript}
\usepackage{yfonts}
\usepackage{tikz}
\usetikzlibrary{arrows.meta}
\usetikzlibrary{shapes,calc}
\usepackage{makeidx}
\usepackage{stmaryrd}
\usepackage{multicol}
\usepackage{bm}
\usepackage[all]{xy}
\numberwithin{equation}{section}
\newtheorem{thm}{Theorem}[section]
\newtheorem{prop}[thm]{Proposition}
\newtheorem{lem}[thm]{Lemma}
\newtheorem{cor}[thm]{Corollary}
{\bf}{\it}

\newtheorem{fthm}{Theorem}{\bf}{\it}
{\bf}{\it}
\newtheorem{fcor}[fthm]{Corollary}{\bf}{\it}
{\bf}{\it}
{\bf}{\it}

\theoremstyle{definition}
\newtheorem{defn}[thm]{Definition}

{\bf}{\rm}

\theoremstyle{remark}
\newtheorem{ex}[thm]{Example}
\newtheorem{rem}[thm]{Remark}
{\bf}{\it}

\newtheorem{definition and corollary}[thm]{Definition and Corollary}

{\it}{\rm}

\newcommand{\A}{{\mathbb A}}
\DeclareMathOperator{\gl}{\mathfrak{gl}}
\DeclareMathOperator{\GL}{GL}
\DeclareMathOperator{\SL}{SL}

\DeclareMathOperator{\Hom}{Hom}
\newcommand{\bB}{{\mathbf B}}
\newcommand{\tbB}{\widetilde{\mathbf B}}
\newcommand{\tbG}{\widetilde{\mathbf G}}

\newcommand{\af}{\mathrm{af}}
\newcommand{\bP}{{\mathbf P}}
\newcommand{\tbP}{\widetilde{\mathbf P}}
\newcommand{\C}{{\mathbb C}}

\newcommand{\cE}{{\mathcal E}}

\newcommand{\tL}{{\mathtt L}}
\newcommand{\cO}{{\mathcal O}}
\newcommand{\bG}{{\mathbf G}}
\newcommand{\Par}{{\mathtt{Par}}}

\newcommand{\bi}{{\mathbf i}}

\DeclareMathOperator*{\gch}{\mathrm{gch}}

\DeclareMathOperator{\Lie}{\mathrm{Lie}}
\DeclareMathOperator{\Pic}{\mathrm{Pic}}
\newcommand{\ro}{\mathrm{rot}}

\newcommand{\la}{\lambda}
\newcommand{\La}{\Lambda}

\newcommand{\Spec}{\mathrm{Spec}}

\newcommand{\Sym}{\mathfrak{S}}
\newcommand{\tSym}{\widetilde{\mathfrak{S}}}

\newcommand{\g}{\mathfrak{g}}
\newcommand{\tg}{\widetilde{\mathfrak{g}}}

\newcommand{\tb}{\widetilde{\mathfrak{b}}}

\newcommand{\gs}{\mathfrak{s}}

\newcommand{\wgt}{\widetilde{\mathfrak{t}}}
\newcommand{\tT}{\widetilde{T}}

\newcommand{\tI}{\mathtt{I}}
\newcommand{\te}{\mathtt{e}}
\newcommand{\ti}{\mathtt{i}}
\newcommand{\hs}{\mathtt{h}}

\newcommand{\sP}{\mathsf{P}}
\newcommand{\sfQ}{\mathsf{Q}}

\newcommand{\gn}{\mathfrak{n}}

\newcommand{\gp}{\mathfrak{p}}
\newcommand{\tp}{\widetilde{\mathfrak{p}}}
\newcommand{\gu}{\mathfrak{u}}
\newcommand{\bv}{\mathbf{v}}

\renewcommand{\P}{\mathbb{P}}

\newcommand{\Proj}{\mathrm{Proj}}
\newcommand{\bQ}{\mathbf{Q}}

\newcommand{\sX}{\mathscr{X}}
\newcommand{\sY}{\mathscr{Y}}

\newcommand{\Z}{\mathbb{Z}}

\newcommand{\Span}{\mbox{\rm Span}}

\newcommand{\Ga}{\mathbb G_a}
\newcommand{\Gm}{\mathbb G_m}

\title{A geometric realization of Catalan functions}
\author{Syu \textsc{Kato}\footnote{Department of Mathematics, Kyoto University, Oiwake Kita-Shirakawa Sakyo Kyoto 606-8502 JAPAN \texttt{E-mail:syuchan@math.kyoto-u.ac.jp}}}
\date{\today}
\begin{document}
\maketitle

\begin{abstract}
We construct smooth projective varieties $\sX_\Psi$, each of which compactifies an equivariant vector subbundle of the cotangent bundle of the flag variety for $\GL(n)$. Each variety $\sX_\Psi$ carries a natural family of line bundles whose spaces of global sections realize, as graded characters, the Catalan functions introduced in Chen's thesis and further studied by Blasiak, Morse, Pun, and Summers. Using the geometry of $\sX_\Psi$, we prove the Chen--Haiman vanishing conjecture and confirm the tame case of the Blasiak--Morse--Pun vanishing conjecture. We further establish the Shimozono--Weyman monotonicity conjectures.
\end{abstract}

\section*{Introduction}

In 2003, Lapointe, Lascoux, and Morse~\cite{LLM03} introduced the notion of $k$-Schur functions in pursuit of a deeper understanding of the internal structure of Macdonald polynomials~\cite{Mac95,Hai01}. These functions were subsequently shown to represent Schubert classes in the affine Grassmannian~\cite{Lam06}, and thus arise naturally in the study of the quantum cohomology of the flag variety $X$ associated with $G = \GL(n, \C)$~\cite{Pet97, LS10}. Nevertheless, their precise relation to Macdonald polynomials, as well as their role in explicit computations in quantum cohomology, remains only partially understood.

In this context, Chen~\cite{Che10} formulated a series of striking conjectures concerning the internal structure of $k$-Schur functions and their generalizations, known as Catalan functions. These conjectures arise from a geometric framework involving certain equivariant vector bundles on the flag variety~$X$. As special cases, they include a conjectural answer to a problem of Broer~\cite[3.16]{Bro94} in type~$\mathsf{A}$, as well as the Shimozono--Weyman vanishing conjecture~\cite{SW00}. While the corresponding character formulas and positivity results have been established by Blasiak, Morse, Pun, and Summers~\cite{BMPS,BMP}, the cohomological components, formulated as the Chen--Haiman vanishing conjecture and its extensions, remain open. These conjectures form the foundation of the geometric program initiated in~\cite{Che10}, and likewise underlie the structure of the monotonicity conjectures in~\cite[\S2.10]{SW00}. In this light, the vanishing results can be seen as the final pieces in a conceptual framework that has taken shape over decades of work by Chen and Haiman, Shimozono and Weyman, and others.

In this paper, we define and study a smooth projective variety $\sX_\Psi$ that compactifies the $G$-equivariant vector subbundle $T^*_\Psi X \subset T^*X$ introduced in~\cite{Che10}.

To state our results precisely, we begin by fixing notation. Let $\Psi$ be a Dyck path of size $n$, corresponding to a root ideal of type~$\mathsf{A}_{n-1}$~\cite{Cel00}, which specifies the subbundle $T^*_\Psi X$. Let $\Par$ denote the set of partitions of length at most $n$, which parametrizes the irreducible polynomial representations of $G$ up to isomorphism. For each $\la \in \Par$, let $V(\la)$ denote the corresponding representation, whose character is the Schur polynomial $s_\la$. Encoding the $\C^\times$-weights as powers of $q$, we write $\gch V$ for the graded character of a rational $(G \times \C^\times)$-module $V$. For any such module $M$, we denote by $M^\vee$ its restricted dual, namely the direct sum of the duals of its $\C^\times$-isotypic components.

The Catalan symmetric function associated with a Dyck path $\Psi$ of size $n$ and $\la \in \Par$ is defined as the Euler-Poincar\'e characteristic
\begin{align}\label{eqn:origCat}
H(\Psi; \la) := {} & \\
\sum_{i \in \Z, \mu \in \Par, m \in \Z} \! & \, (-1)^i q^m s_\mu \cdot \dim \Hom_{G \times \C^\times} \bigl( V(\mu) \boxtimes \C_{-m\delta}, H^i\bigl(T^*_\Psi X, \cO_{T^*_\Psi X}(\la)\bigr)^\vee \bigr), \nonumber
\end{align}
where $H(\Psi; \la) = H(\Psi; \la; w_0)$ in~\cite[(2.2)]{BMP}, and $\C_{m\delta}$ denotes the one-dimensional $\C^\times$-representation of weight~$m$. The Chen--Haiman vanishing conjecture asserts that the higher cohomology groups vanish, in which case~\eqref{eqn:origCat} coincides with
\begin{equation}
\sum_{\mu \in \Par, m \in \Z}  q^m s_\mu \cdot \dim \Hom_{G \times \C^\times} \bigl( V(\mu) \boxtimes \C_{-m\delta}, H^0\bigl(T^*_\Psi X, \cO_{T^*_\Psi X}(\la)\bigr)^\vee \bigr).\label{eqn:H0Cat}
\end{equation}
The sum in~\eqref{eqn:H0Cat} is finite, although the ambient space satisfies
\[
\dim\, H^0\bigl(T^*_\Psi X, \cO_{T^*_\Psi X}(\la)\bigr) = \infty
\]
in general. Most of the irreducible rational representations of $G$ occurring in $H^0(T^*_\Psi X, \cO_{T^*_\Psi X}(\la))$ are therefore not captured by~\eqref{eqn:H0Cat}; these missing components correspond precisely to the rational but non-polynomial representations of $G$.

Our main results are summarized below.

\begin{fthm}[$\doteq$ Theorems~\ref{thm:str},
\ref{thm:XTid}, \ref{thm:Xmain}, and \ref{thm:incl}]\label{fthm:main}
For every root ideal $\Psi\subset\Delta^+$, the variety $\sX_\Psi$ is a
smooth projective $(G\times\C^\times)$-variety with the following properties.
\begin{enumerate}
\item There exists a $(G \times \C^\times)$-equivariant open immersion
$
T^*_\Psi X \hookrightarrow \sX_\Psi.
$
\item For each $\la \in \Par$, there exists a $(G \times \C^\times)$-equivariant line bundle $\cO_{\sX_\Psi}(\la)$ on $\sX_\Psi$ such that
\begin{align*}
H^{>0}(\sX_\Psi, \cO_{\sX_\Psi}(\la)) &= 0, \\
\gch H^0(\sX_\Psi, \cO_{\sX_\Psi}(\la))^\vee &= \left[ H(\Psi; \la) \right]_{q \mapsto q^{-1}}.
\end{align*}
\item There is an effective $(G\times\C^\times)$-equivariant Cartier divisor
$\partial$, with support $\sX_\Psi\setminus T^*_\Psi X$, such that
\[
 H^{>0}\bigl(
   \sX_\Psi,\cO_{\sX_\Psi}(\la+m\partial)
 \bigr)=0
 \qquad(\la\in\Par,\ m\ge0).
\]
Consequently,
\begin{equation}
 H^{>0}\bigl(T^*_\Psi X,\cO_{T^*_\Psi X}(\la)\bigr)=0
 \qquad(\la\in\Par).\label{eqn:fTvan}
\end{equation}
A parabolic analogue of this vanishing result also holds; see Corollary~\ref{cor:pv}.
\end{enumerate}
\end{fthm}

Theorem~\ref{fthm:main}(3) resolves the Chen--Haiman vanishing conjecture~\cite[Conjecture 5.4.3(2)]{Che10}. Combined with~\cite[Theorem~2.18]{BMP}, this establishes~\cite[Conjecture~5.4.3]{Che10} in full generality. Since this conjecture provides an answer to a question of Broer~\cite[3.16]{Bro94} (in type~$\mathsf{A}$) and extends the Shimozono--Weyman vanishing conjecture~\cite[\S2.4]{SW00}, our result settles both as well (see Remark~\ref{rem:coverage}). In the special case $T^*_\Psi X = T^*X$, the variety $\sX_\Psi$ recovers the smooth resolution~\cite{Ngo99, MV03} of Lusztig's compactification~\cite{Lus81b} of the nilpotent cone of~$\mathfrak{gl}(n, \C)$. For completeness, we also note in Remark~\ref{rem:gen} that our argument applies over fields of positive characteristic, with suitable modifications.

The construction of $\sX_\Psi$ starts from the Blasiak--Morse--Pun rotation formula~\cite[Theorem~2.3]{BMP}. Lifting its Demazure-operator expression to Demazure functors realizes the restricted duals of the resulting modules as spaces of sections of line bundles on a suitable Bott--Samelson variety. Taking the multigraded Proj of the resulting section algebra gives $\sX_\Psi$. A nontrivial reorganization of the Demazure functors then reveals an iterated projective-space bundle structure on $\sX_\Psi$, which is essential for proving Theorem~\ref{fthm:main}.

As a corollary of Theorem~\ref{fthm:main}, we obtain the following:

\begin{fcor}[$\doteq$ Lemma~\ref{lem:infiinj}]\label{fcor:inj}
There exists an action of $\GL(n, \C[\![z]\!]) \rtimes \C^\times$ on $\sX_\Psi$ such that the natural restriction map
\[
H^0(\sX_\Psi, \cO_{\sX_\Psi}(\la)) \hookrightarrow H^0(T^*_{\Psi} X, \cO_{T^*_{\Psi} X}(\la)), \qquad \la \in \Par,
\]
is an inclusion of graded representations of $\mathfrak{gl}(n, \C[z])$.
\end{fcor}

A local chart analysis of $\sX_\Psi$ further yields the following result:

\begin{fthm}[$\doteq$ Theorem~\ref{thm:infi}]\label{fthm:head}
For each $\la \in \Par$, the space $H^0(\sX_\Psi, \cO_{\sX_\Psi}(\la))$ has a simple head $V(\la)^{\vee}$ as a graded $\mathfrak{gl}(n, \C[z])$-module.
\end{fthm}

As an additional consequence of our construction, we obtain the following:

\begin{fcor}[$\doteq$ Corollary~\ref{cor:SW}]\label{fcor:SW}
Let $\Psi' \subset \Psi$ be an inclusion of Dyck paths, which in particular induces an inclusion $T^*_{\Psi'} X \subset T^*_\Psi X$. Then for each $\la \in \Par$, the restriction map
\[
H^0(T^*_\Psi X, \cO_{T^*_\Psi X}(\la)) \longrightarrow H^0(T^*_{\Psi'} X, \cO_{T^*_{\Psi'} X}(\la))
\]
is surjective.
\end{fcor}

In \S\ref{subsec:mon}, we reformulate~\cite[Conjectures~12 and~13]{SW00}, along with their natural generalizations, as module-theoretic statements, and establish them using Corollary~\ref{fcor:SW}.

In the statements above and throughout the introduction,
$\sX_\Psi$ denotes $\sX_\Psi(w_0)$, where
$w_0\in\Sym_n$ is the longest element. We have a subvariety $\sX_\Psi ( w ) \subset \sX_\Psi$ for each $w \in \Sym_n$. From this, the vanishing result~\eqref{eqn:fTvan} extends to tame $w \in \Sym_n$ (\cite{BMP}), which proves the tame case of~\cite[Conjecture~3.4(ii)]{BMP}. In a similar manner, Corollary~\ref{cor:BMP} proves the excellent-filtration statement in~\cite[Conjecture~3.4(iv)]{BMP}. Note that we can replace an arbitrary $w$ with tame $v  \in \Sym_n$ such that $\sX_\Psi ( w ) = \sX_\Psi ( v )$ (essentially by Lemma~\ref{lem:tame}). Hence we need to analyze a suitable subvariety of $\sX_\Psi (w)$ in order to establish~\cite[Conjecture~3.4(ii)(iii)]{BMP} in its full generality. The corresponding vanishing claim for arbitrary $w\in\Sym_n$ made in an
earlier version of this paper was incorrect and is withdrawn here.

The organization of this paper is as follows.
In Section~1, we fix notation and recall the necessary facts about
root ideals, affine Demazure functors and modules, and affine flag
varieties.
Section~2 is devoted to a new expression of the rotation theorem from~\cite{BMP}.
In Section~3, we construct the variety $\sX_\Psi$ (Theorem~\ref{thm:str}) and illustrate it with an explicit example (Example~\ref{ex:n4}).
In Section~4, we establish Theorem~\ref{fthm:main}(1)(2).
Section~5 proves Theorem~\ref{fthm:main}(3) and derives the
infinitesimal-action, simple-head, and monotonicity results stated above.

The varieties introduced here serve as natural geometric counterparts to the Catalan functions. A promising direction for future research is to place these constructions within the framework of topological field theories and the geometric realization of Macdonald polynomials associated with $G = \GL(n)$. We hope to return to these questions in future work.

\section{Preliminaries}\label{sec:prelim}

Throughout the main body of this paper, we work over the field $\C$ of complex numbers. A \emph{variety} means a separated, integral, normal scheme of finite type over~$\C$. When the topology and scheme structure are clear from context, we often identify a variety $\sX$ with its set of $\C$-points, denoted $\sX(\C)$. In particular, we write $\Gm$ and $\Ga$ for the multiplicative group~$\C^\times$ and the additive group~$\C$, respectively.

Let $\tL$ be a free abelian monoid. A $\tL$-graded vector space is a $\C$-vector space $V$ equipped with a direct sum decomposition $V = \bigoplus_{a \in \tL} V_a$, where each graded piece $V_a$ is finite-dimensional. Given such a grading, we define the restricted dual as $V^{\vee} := \bigoplus_{a \in \tL} V_a^*$. A $\tL$-graded ring is a unital $\C$-algebra $R = \bigoplus_{a \in \tL} R_a$ satisfying $\C \cdot 1 = R_0$ and $R_a \cdot R_{a'} \subset R_{a + a'}$ for all $a, a' \in \tL$.

Suppose that $R$ is commutative. The grading defines an action of the torus
\[
T_\tL:=\Hom_{\mathrm{mon}}(\tL,\Gm)\cong\Gm^r
\]
on $\Spec R$, where $\tL\cong\Z_{\ge0}^r$. Let
$\mathrm{Irr}(R) \subset \Spec R$ be the complement of the locus of points with finite stabilizer. We use the multigraded version of $\Proj$ defined as
\begin{equation}\label{eqn:mproj}
\Proj_\tL R := \left( \Spec R \setminus \mathrm{Irr}(R) \right) / T_\tL,
\end{equation}
where the right-hand side is a geometric quotient.
Equivalently, this is the union of the opens determined by the
relevant homogeneous elements; see~\cite[\S2]{BS03}.
We will apply this construction only to the $\Par$-graded rings in
Section~\ref{sec:XP}. In that setting, the action on the relevant
locus is free, and the base-point-free linear systems constructed
there realize the quotient inside a product of projective spaces.

We also use the relative version $\Proj_{S,\tL}\mathcal A$ for a
$\tL$-graded sheaf of algebras on a scheme $S$; it is defined
locally by the above construction and is compatible with base change.

Given a representation~$M$, we define its \emph{head} to be the largest semisimple quotient of~$M$.

For general background, we refer the reader to the standard references~\cite{Kum02,CG97}.

\subsection{Algebraic groups}

We fix an integer $n > 1$ and define the algebraic group
\[
G := \C^{\times} \mathrm{Id} \cdot \SL(n) = \GL(n) \subset M_n \cong \C^{n^2}.
\]
We also define the proalgebraic group $\bG := \C^{\times} \mathrm{Id} \cdot \SL(n, \C[\![z]\!])$ over~$\C$. In addition, we consider the group
\[
G(\!(z)\!) := \C^{\times} \mathrm{Id} \cdot \SL(n,\C(\!(z)\!)),
\]
regarded as a topological group.

Let $E_{ij} \in M_n$ ($1 \le i, j \le n$) denote the standard matrix units. Let $T \subset G$ be the diagonal torus, and let $B \subset G$ (resp.\ $B^- \subset G$) denote the subgroup of upper (resp. lower) triangular matrices in $G$. The unipotent subgroup $N := [B, B] \subset B$ consists of upper unitriangular matrices in $G$.

We have the evaluation map
\[
\mathtt{ev}_0 : \bG \longrightarrow G \qquad f( z ) \mapsto f( 0 ).
\]
We define $\bB := \mathtt{ev}_0^{-1}(B)$.

For each $1 \le i < n$, let $P_i \subset G$ denote the (algebraic) subgroup generated by $B$ and $\mathrm{Id} + \C E_{i+1,i}$, and let $\bP_i \subset \bG$ denote the proalgebraic subgroup generated by $\bB$ and $\mathrm{Id} + \C E_{i+1,i}$. We define $\bP_0$ to be the (pro)algebraic subgroup of $G(\!(z)\!)$ generated by $\bB$ and $\mathrm{Id} + \C E_{1,n} z^{-1}$.

Observe that there is a loop rotation $\Gm$-action (denoted $\Gm^\ro$) on each of $\bB$, $\bP_i$, and $\bG$.

We denote by $\widehat{\bB}$, $\widehat{\bP}_i$, and $\widehat{\bG}$ the semidirect products of $\bB$, $\bP_i$, and $\bG$ with $\Gm^\ro$, respectively. In addition, the group $G(\!(z)\!)$ admits a central extension by $\C^\times$, which induces a trivial central extension $\widetilde{\bP}_i$ ($0 \le i < n$) of $\widehat{\bP}_i$ by $\Gm$. We denote this copy of $\Gm$ by $\Gm^{\mathrm{ce}}$.

We define the extended torus
\[
\widehat{T} := T \times \Gm^\ro \times \{ 1 \} \subset T \times \Gm^\ro \times \Gm^{\mathrm{ce}} =: \tT,
\]
so that $\tbB := \widehat{\bB} \times \Gm^{\mathrm{ce}}$ contains $\tT$, and $\widehat{\bB} \cap \tT = \widehat{T}$. We also set
\[
\tbG := \widehat{\bG} \times \Gm^{\mathrm{ce}} \supset \tbB, \widehat{\bG}, \quad \text{with} \quad \tbB \cap \widehat{\bG} = \widehat{\bB}.
\]
Moreover, we have $\tbP_i \cap \tbP_j = \tbB$ whenever $i \ne j$. For each $0 \le i < n$, there exists a unique $\tT$-stable algebraic subgroup of $\tbP_i$ isomorphic to $\SL(2)$, which we denote by $\SL(2, i)$.

We denote the Lie algebra of an algebraic group by the corresponding lowercase German letter.

We define the group
\[
\widetilde{G}(\!(z)\!) := \Gm^\ro \ltimes G(\!(z)\!) \ltimes \Gm^{\mathrm{ce}}.
\]
Let $\tbG^- \subset \widetilde{G}(\!(z)\!)$ denote the subgroup generated by $\tT \cdot G$ and $\mathrm{Id} + \C E_{1,n} z^{-1}$.

Note that the groups $\widetilde{G}(\!(z)\!)$ and $\tbG^-$ are not algebraic.

\subsection{Weights and Weyl group actions}

For each $1 \le i \le n$, we define the algebraic character $\epsilon_i : T \to \Gm$ that extracts the $i$-th (diagonal) entry of an element of $T$. We set $\sP := \bigoplus_{i=1}^n \Z \epsilon_i$. Consider the subsets of polynomial and dominant weights given by
\[
\mathtt{Comp} := \sum_{i=1}^n \Z_{\ge 0} \epsilon_i, \qquad \sP^+ := \left\{ \sum_{i=1}^n \la_i \epsilon_i \in \sP \middle| \la_1 \ge \la_2 \ge \cdots \ge \la_n \right\}.
\]
For a weight $\la = \sum_{i=1}^n \la_i \epsilon_i \in \sP$, we define its total weight as $|\la| := \sum_{i=1}^n \la_i \in \Z$.

The symmetric group $\Sym_n$ acts on $\sP$ and $\mathtt{Comp}$ by permuting the indices.

We define $\Par := ( \sP^+ \cap \mathtt{Comp} )$ and identify it with the set of partitions of length at most $n$. The semigroup $\Par$ is generated by the elements
$$\varpi_i := \epsilon_1 + \cdots + \epsilon_i \hskip 5mm 1 \le i \le n.$$

Let $\wp$ and $\delta$ denote the degree-one characters of $\Gm^{\mathrm{ce}}$ and $\Gm^\ro$, respectively, each extended trivially to $\tT$. We regard $\varpi_i$ as a character of $\tT$ via the projection $\tT \to T$, and refer to this as the \emph{standard lift} of $\varpi_i$.

We define an alternative, \emph{nonstandard lift} of $\varpi_i$ to $\tT$ by setting
\begin{equation}\label{eqn:nonst}
\La_i := \varpi_i + \wp \qquad (1 \le i \le n).
\end{equation}
This extends linearly to define a nonstandard lift of any character of $T$ to $\tT$.

We set $\tI_\af:=\{0,1,\ldots,n-1\}$ and
$\tI:=\{1,2,\ldots,n-1\}$. These are the index sets of the affine Dynkin diagram of type $\mathsf{A}_{n-1}^{(1)}$ and the Dynkin diagram of type $\mathsf{A}_{n-1}$, respectively.
When an affine Dynkin index labels one of the families
$\{\varpi_i\}$ or $\{\La_i\}$, we represent $0\in\tI_\af$ by the
subscript $n$. This is only an indexing convention: the level-one
fundamental weight at node $0$ is $\wp$, whereas
$\La_n=\varpi_n+\wp$. Note that the sets $\{\varpi_i\}_{i \in \tI_\af}$ and $\{\La_i\}_{i \in \tI_\af}$ correspond to each other via restriction to $T$.

We define the affine weight lattice and its subset of dominant weights by
\[
\sP_\af := \bigoplus_{i=1}^n \Z \varpi_i \oplus \Z \wp \oplus \Z \delta,
\qquad
\sP^+_\af := \left( \sum_{i =1}^{n-1} \Z_{\ge 0} \La_i \right) + \Z_{\ge 0} \wp + \Z \varpi_n + \Z \delta \subset \sP_\af.
\]
We identify $\sP_\af$ with the character lattice of $\tT$.

The set of positive roots of $G$ is given by $\Delta^+ := \{ \epsilon_i - \epsilon_j \mid 1 \le i < j \le n \} \subset \sP$. We define $\alpha_i := \epsilon_i - \epsilon_{i+1}$ for $1 \le i < n$, and set $\alpha_0 := -\vartheta + \delta$, where $\vartheta := \epsilon_1 - \epsilon_n$ denotes the highest root.

We equip $\sP_\af$ with a symmetric bilinear form defined by
\[
\langle \epsilon_i, \epsilon_j \rangle = \delta_{ij}, \qquad \wp, \delta \in \mathrm{Rad} \, \langle \cdot, \cdot \rangle.
\]

Let $\mathfrak{n} := \Lie N \subset M_n$. For each positive root $\alpha = \epsilon_i - \epsilon_j \in \Delta^+$, we set
\[
\mathfrak{g}_\alpha := \C E_{ij} \subset \mathfrak{n} \subset M_n.
\]

The root lattice $\sfQ \subset \sP$ is defined by $\sfQ := \sum_{\beta \in \Delta^+} \Z \beta$. The permutation action of $\Sym_n$ on $\sP$ restricts to an action on $\sfQ$, and we define
\[
\tSym_n := \Sym_n \ltimes \sfQ.
\]
The standard embedding $\Sym_n \subset G$ via permutation matrices extends naturally to an embedding $\tSym_n \hookrightarrow G(\!(z)\!)$, under which an element $\sum_{i=1}^n \mu_i \epsilon_i \in \sfQ$ is mapped to
\[
z^\mu := 
\begin{pmatrix}
z^{\mu_1} & 0 & \cdots & 0 \\
0 & z^{\mu_2} & \cdots & 0 \\
\vdots & \vdots & \ddots & \vdots \\
0 & 0 & \cdots & z^{\mu_n}
\end{pmatrix}
\in G(\!(z)\!), \qquad \text{where } \sum_{i=1}^n \mu_i = 0.
\]

The group $\tSym_n$ is generated by the elements $\{s_i\}_{i \in \tI_\af}$, where
\[
s_i := 
\begin{cases}
(i, i+1) & \text{if } 1 \le i < n, \\
(1,n) \cdot z^{-\vartheta} & \text{if } i = 0.
\end{cases}
\]
Each $s_i$ lies in the subgroup $\tbP_i$. The group $\tSym_n$ acts on $\sP_\af$ via the rule
\[
s_i(\La) := \La - \left( \langle \alpha_i, \La \rangle + \delta_{i0} \La(K) \right) \alpha_i,
\qquad i \in \tI_\af,
\]
where $K \in \Hom(\sP_\af, \Z)$ is defined by
\[
\varpi_i(K) = 0 \quad (i \in \tI_\af), \qquad \delta(K) = 0, \qquad \wp(K) = 1.
\]

Elements in the $\tSym_n$-orbit of $\{\alpha_i\}_{i \in \tI_\af} \subset \sP_\af$ are called \emph{affine roots}. An affine root is said to be positive if it lies in the semigroup $\sum_i \Z_{\ge 0} \alpha_i$.

Note that the Dynkin diagram automorphism of type $\mathsf{A}_{n-1}^{(1)}$ acts on the set of affine roots (and on the subset of positive affine roots) via the linear transformation that cyclically shifts indices modulo~$n$. This induces an automorphism of $\widetilde{G}(\!(z)\!)$ that fixes scalar matrices.

Every element $w \in \tSym_n$ admits an expression of the form
\begin{equation}\label{eqn:reduced}
w = s_{i_1} s_{i_2} \cdots s_{i_\ell}, \qquad i_1, \dots, i_\ell \in \tI_\af.
\end{equation}
Let $\bi := (i_1, i_2, \ldots, i_\ell)$ denote the sequence of indices appearing in~\eqref{eqn:reduced}. If the length $\ell$ of $\bi$ is minimal among all such expressions for $w$, we call $\bi$ a \emph{reduced expression} of $w$, and refer to $\ell$ as the \emph{length} of $w$.

We define the (strong) Bruhat order on $\tSym_n$ by declaring that $w < v$ if some reduced expression of $w$ appears as an ordered subword of a reduced expression of $v$. The length of an element $w \in \tSym_n$ is denoted by $\ell(w)$.

Let $w_0 \in \Sym_n$ denote the longest element, defined by $w_0(i) = n - i + 1$ for $1 \le i \le n$.

\subsection{Root ideals}

\begin{defn}[Root ideals]\label{def:ri}
A subset $\Psi \subset \Delta^+$ is called a \emph{root ideal} if and only if
\[
(\Psi + \Delta^+) \cap \Delta^+ \subset \Psi.
\]
Equivalently, $\Psi$ is a root ideal if, for every $\epsilon_i - \epsilon_j \in \Psi$, the elements $\epsilon_{i'} - \epsilon_j$ and $\epsilon_i - \epsilon_{j'}$ also lie in $\Psi$ for all $i' < i$ and $j < j'$.

Given a root ideal $\Psi \subset \Delta^+$, we define
\[
\mathfrak{n}(\Psi) := \bigoplus_{\alpha \in \Psi} \mathfrak{g}_\alpha \subset \mathfrak{n}.
\]
We denote by $|\Psi|$ the cardinality of $\Psi$, which coincides with $\dim \mathfrak{n}(\Psi)$.
\end{defn}

\begin{defn}\label{def:rideal}
Let $\Psi \subset \Delta^+$ be a root ideal. For each $1 \le k \le n$, define
\begin{equation}\label{eqn:hk}
\hs_k(\Psi)
:= \#\{1\le i<k\mid E_{ik}\in\mathfrak n(\Psi)\}.
\end{equation}
Thus, $\hs_k(\Psi)$ is the number of roots in $\Psi$ of the form
$\epsilon_i-\epsilon_k$. Equivalently, when $\Delta^+$ is represented
by the staircase array in which the $(i,j)$-box corresponds to
$\epsilon_i-\epsilon_j$, the number $\hs_k(\Psi)$ is the number of
boxes belonging to $\Psi$ in the $k$-th column.

For each $1 \le i \le n$, define
\[
d_i(\Psi) := \# \{ i \le j \le n \mid E_{ij} \not\in \mathfrak{n}(\Psi) \},
\qquad
e_i(\Psi) := i + d_i(\Psi).
\]
We set
\[
\tI(\Psi) := \{ 1 \le i < n \mid e_i(\Psi) \le n,\ d_i(\Psi) \le d_{i+1}(\Psi) \},
\qquad
\ell(\Psi) := |\tI(\Psi)|.
\]

Let $\{ e_i(\Psi) \mid i \in \tI(\Psi) \}$ be the set of values $e_i(\Psi)$ indexed by $\tI(\Psi)$, and let
\[
\{ e_i(\Psi) \}_{i \in \tI(\Psi)} = \{ \mathtt{e}_1(\Psi) < \mathtt{e}_2(\Psi) < \cdots < \mathtt{e}_\ell(\Psi) \}, \quad \mathtt{e}_{\ell+1}(\Psi) := e_{n+1}(\Psi) := n+1,
\]
denote its increasing rearrangement. For each $1 \le j \le \ell(\Psi)$, there exists a unique $i \in \tI(\Psi)$ such that $\mathtt{e}_j(\Psi) = e_i(\Psi)$, and we define $\mathtt{i}_j(\Psi) := i$. By convention, we also set $\mathtt{i}_0(\Psi) := 0$, and $e_0(\Psi) = \te_0(\Psi) := 1$.
\end{defn}

The invariants in Definition~\ref{def:rideal} are illustrated in
Example~\ref{ex:n6}.

\begin{defn}[$\Psi$-tame elements]\label{def:tame}
Let $\Psi \subset \Delta^+$ be a root ideal. An element $w \in \Sym_n$ is said to be \emph{$\Psi$-tame} if $w s_i < w$ for all $i$ with $d_1(\Psi) < i < n$.

We define $w_0^\Psi$ to be the longest element in the subgroup
\[
\Sym_{n - d_1(\Psi)} \cong \langle s_{e_1(\Psi)}, s_{e_1(\Psi)+1}, \ldots, s_{n-1} \rangle \subset \Sym_n.
\]
\end{defn}

\begin{ex}\label{ex:n6}
Assume that $n = 6$, and consider the root ideal
\[
\Psi = \{ \epsilon_1 - \epsilon_3,\, \epsilon_1 - \epsilon_4,\, \epsilon_1 - \epsilon_5,\, \epsilon_1 - \epsilon_6,\,
\epsilon_2 - \epsilon_3,\, \epsilon_2 - \epsilon_4,\, \epsilon_2 - \epsilon_5,\, \epsilon_2 - \epsilon_6,\,
\epsilon_3 - \epsilon_6 \}.
\]
We compute:
\[
d_1(\Psi) = 2,\quad d_2(\Psi) = 1,\quad d_3(\Psi) = 3,\quad d_4(\Psi) = 3,\quad d_5(\Psi) = 2,\quad d_6(\Psi) = 1,
\]
and hence $\quad e_1(\Psi) = 3,\quad e_2(\Psi) = 3,\quad e_3(\Psi) = 6,\quad e_4(\Psi) = 7$. This implies
\[
\te_1(\Psi) = 3,\quad \te_2(\Psi) = 6,\quad \ell(\Psi) = 2.
\]

Since $e_5(\Psi) = e_6(\Psi) = 7 > n = 6$, these values do not contribute to the sets $\{ \te_j(\Psi) \}$, $\{ \ti_j(\Psi) \}$, or to $\ell(\Psi)$.  
By convention, we have $\ti_0(\Psi) = 0$, and from the values above, we find $\ti_1(\Psi) = 2$ and $\ti_2(\Psi) = 3$.  
Hence, $\tI(\Psi) = \{ \ti_1(\Psi), \ti_2(\Psi) \} = \{ 2, 3 \}$.

For $3 = e_1(\Psi) \le k \le n$, we compute:
\[
\hs_3(\Psi) =
\hs_4(\Psi) =
\hs_5(\Psi) = \ti_1(\Psi) = 2, \quad
\hs_6(\Psi) = \ti_2(\Psi) = 3,
\]
using that $\te_1(\Psi) = e_2(\Psi) = 3$ and $\te_2(\Psi) = e_3(\Psi) = 6$.

This situation is illustrated in the diagram below:
\begin{center}
\begin{tikzpicture}[scale=0.7, thick, >=Stealth]

  \def\s{1.0}

  \draw[thick] (0,0) rectangle (6*\s,6*\s);

  \foreach \i in {0,...,5} {
    \pgfmathsetmacro{\x}{6 - \i - 1}
    \pgfmathsetmacro{\y}{\i}
    \draw (\x*\s,\y*\s) rectangle ++(\s,\s);
  }

  \foreach \x/\y in {
    2/4, 2/5,
    3/4, 4/4, 5/4,
    3/5, 4/5, 5/5,
    5/3} {
    \fill[red!30] (\x*\s,\y*\s) rectangle ++(\s,\s);
    \draw[thick, dash pattern=on 2pt off 5pt] (\x*\s,\y*\s) rectangle ++(\s,\s);
  }

  \node at (4.5*\s,4.6*\s) {$\Psi$};

  \node[left] at (-0.2,6*\s) {$\mathtt{i}_0$};
  \node[left] at (-0.2,4*\s) {$\mathtt{i}_1$};
  \node[left] at (-0.2,3*\s) {$\mathtt{i}_2$};

  \draw[<->, shorten >=5pt, shorten <=5pt]
    (0,4.5*\s) -- (3*\s,4.5*\s)
    node[midway, fill=white, inner sep=1pt] {$\mathtt{e}_1$};

  \draw[<->, shorten >=5pt, shorten <=5pt]
    (2*\s,3.5*\s) -- (5*\s,3.5*\s)
    node[midway, fill=white, inner sep=1pt] {$d_3$};

  \draw[<->, shorten >=5pt, shorten <=5pt]
    (0*\s,5.5*\s) -- (2*\s,5.5*\s)
    node[midway, fill=white, inner sep=1pt] {$d_1$};

  \draw[<->, shorten >=5pt, shorten <=5pt]
    (0,1.5*\s) -- (7*\s,1.5*\s)
    node[midway, fill=white, inner sep=1pt] {$e_5$};
    
  \draw[<->, shorten >=5pt, shorten <=5pt]
    (5.5*\s,3*\s) -- (5.5*\s,6*\s)
    node[midway, fill=white, inner sep=1pt] {$\mathtt{h}_6$};

  \draw[<->, shorten >=5pt, shorten <=5pt]
    (3.5*\s,4*\s) -- (3.5*\s,6*\s)
    node[midway, fill=white, inner sep=1pt] {$\mathtt{h}_4$};

  \draw[thick, dash pattern=on 5pt off 7pt] (0*\s,4*\s) -- (1*\s,4*\s);
  \draw[thick, dash pattern=on 5pt off 7pt] (0*\s,3*\s) -- (2*\s,3*\s);
  \draw[thick, dash pattern=on 5pt off 7pt] (4*\s,3*\s) -- (5*\s,3*\s);
\end{tikzpicture}
\end{center}
The red-shaded boxes represent the elements of $\Psi$.
\end{ex}

We now summarize basic properties of the invariants associated with a root ideal~$\Psi$.

\begin{lem}[Cellini {\cite[\S3]{Cel00}}]\label{lem:rootid}
For any root ideal $\Psi \subset \Delta^+$, the subspace $\mathfrak{n}(\Psi) \subset \mathfrak{n}$ is $B$-stable. Moreover, every $B$-stable subspace of $\mathfrak{n}$ arises uniquely in this way.
\end{lem}

\begin{rem}
As shown in~\cite[\S4]{Pan04}, the set of $B$-stable ideals in $\mathfrak{n}$, that is the set of root ideals, is naturally in bijection with the set of Dyck paths of size~$n$.
\end{rem}

\begin{lem}\label{lem:basicinv}
Let $\Psi \subset \Delta^+$ be a root ideal, and let $1 \le i < n$. Then
\[
d_i ( \Psi ) \le d_{i+1} ( \Psi ) + 1, \qquad \text{and} \qquad i < e_i ( \Psi ) \le e_{i+1} ( \Psi ) \le n+1.
\]
Moreover, for all $1 \le j \le \ell ( \Psi )$, we have
\begin{equation}
\begin{gathered}
\ti_j(\Psi)<\te_j(\Psi),
\qquad
\ti_{j-1}(\Psi)<\ti_j(\Psi),\\
e_k(\Psi)=\te_j(\Psi)
\qquad
\bigl(\ti_{j-1}(\Psi)<k\le\ti_j(\Psi)\bigr).
\end{gathered}
\label{eqn:teblock}
\end{equation}
In addition, we have $e_k (\Psi) = n + 1$ for $\ti_{\ell(\Psi)} (\Psi) < k \le n$.
\end{lem}

\begin{proof}
The inequalities in the first display are straightforward. For
$1\le i<n$, we have
\[
 d_i(\Psi)\le d_{i+1}(\Psi)
 \quad\Longleftrightarrow\quad
 e_i(\Psi)<e_{i+1}(\Psi).
\]
Hence $\tI(\Psi)$ is the set of right endpoints of the maximal
intervals on which $e_\bullet(\Psi)$ is constant with value at most
$n$. The assertions~\eqref{eqn:teblock} follow from the definitions of
$\ti_j(\Psi)$ and $\te_j(\Psi)$.

The final assertion follows from the preceding description of
$\tI(\Psi)$. Indeed, no index
$k>\ti_{\ell(\Psi)}(\Psi)$ satisfies $e_k(\Psi)\le n$.
Since $e_k(\Psi)\le n+1$, it follows that $e_k(\Psi)=n+1$.
\end{proof}

\begin{lem}\label{lem:transp}
Let $\Psi \subset \Delta^+$ be a root ideal. Let $0 \le s \le \ell(\Psi)$, and suppose $\te_s(\Psi) \le j < \te_{s+1}(\Psi)$.
Then
\[
E_{ij} \in \mathfrak{n}(\Psi)
\quad\Longleftrightarrow\quad
1 \le i \le \ti_s(\Psi).
\]
In particular,
\begin{equation}\label{eqn:hk-piecewise}
\hs_j(\Psi)=\ti_s(\Psi)
\qquad
\text{if} \quad \te_s(\Psi)\le j<\te_{s+1}(\Psi).
\end{equation}
\end{lem}

\begin{proof}
If $1 \le i \le \ti_s(\Psi)$, then by Lemma~\ref{lem:basicinv} we have $e_i(\Psi) \le \te_s(\Psi) \le j$, and hence $E_{ij} \in \gn(\Psi)$.
Conversely, if $i > \ti_s(\Psi)$, then Lemma~\ref{lem:basicinv} implies that
$e_i(\Psi)\ge\te_{s+1}(\Psi)>j$. Hence, $E_{ij}\notin\gn(\Psi)$.
The equivalence follows. Counting the admissible indices $i$ gives
\eqref{eqn:hk-piecewise}.
\end{proof}

\begin{lem}\label{lem:count}
Let $\Psi \subset \Delta^+$ be a root ideal, and suppose $1 < k \le n$. Then
\[
\hs_{k-1} ( \Psi ) \le \hs_k ( \Psi ) < k.
\]
Moreover, the cardinality of $\Psi$ is given by
\begin{equation}
|\Psi| = \sum_{k=1}^n \hs_k (\Psi).\label{eqn:Psicard}
\end{equation}
\end{lem}

\begin{proof}
By Definition~\ref{def:rideal}, we have
\[
\hs_k(\Psi)=\#\{i\mid E_{ik}\in\mathfrak n(\Psi)\}.
\]
Since $\Psi$ is a root ideal, $E_{i,k-1}\in\mathfrak n(\Psi)$ implies
$E_{ik}\in\mathfrak n(\Psi)$, and hence
$\hs_{k-1}(\Psi)\le\hs_k(\Psi)$. Moreover,
$E_{ik}\in\mathfrak n(\Psi)\subset\mathfrak n$ implies $i<k$, so
$\hs_k(\Psi)<k$. Finally,~\eqref{eqn:Psicard} follows by summing over
the columns $k=1,\dots,n$.
\end{proof}

\subsection{Representations}\label{subsec:repr}

We call a finite-dimensional rational $G$-module \emph{polynomial} if it is
a polynomial representation of $G=\GL(n)$, equivalently, if all its weights
belong to $\mathtt{Comp}$. Its character then lies in
$\Z[X_1,\ldots,X_n]$, where $X_i=e^{\epsilon_i}$.

For a rational representation $V$ of $\tT$, we define the \emph{graded character} by
\[
\gch V := \sum_{\la \in \sP, \, m \in \Z} q^m e^\la \cdot \dim \Hom_{T \times \Gm^\ro}(\C_{\la + m \delta}, V).
\]

A rational representation of $\tbB$ (resp. $\tbP_i$) is one that factors through a finite-dimensional quotient, and thus defines a rational representation of an algebraic group.

For each $\la \in \sP^+$, let $V(\la)$ denote the irreducible finite-dimensional $G$-module generated by a $B$-eigenvector $\bv_\la$ of $T$-weight $\la$.  
The natural action of $\Sym_n$ on $V(\la)$ then yields a $T$-eigenvector $\bv_{w \la} \in V(\la)$ of weight $w \la \in \sP$ for each $w \in \Sym_n$.

For each $\La \in \sP_\af^+$, let $L(\La)$ denote the integrable highest weight module of $\widetilde{G}(\!(z)\!)$ generated by a $\tbB$-eigenvector $\bv_\La$ of $\tT$-weight $\La$.  
The natural action of $\tSym_n$ on $L(\La)$ then gives rise to a $\tT$-eigenvector $\bv_{w \La} \in L(\La)$ of weight $w \La$ for each $w \in \tSym_n$.

For $\la \in \sP^+$ and $w \in \Sym_n$, the Demazure module of $V(\la)$ is defined by
\[
V_w(\la) := \operatorname{span}_{\C} (  B \bv_{w\la} ) \subset V(\la).
\]
Similarly, for $\La \in \sP_\af^+$ and $w \in \tSym_n$, the Demazure module of $L(\La)$ is defined by
\[
L_w(\La) := \operatorname{span}_{\C} ( \tbB \bv_{w\La} ) \subset L(\La).
\]

\subsection{Geometric interpretation of Demazure functors}\label{subsec:fv}

We set $X := G/B$, the flag variety of $G$.  
For each $\la \in \sP$, let $\cO_X(\la)$ denote the $G$-equivariant line bundle on $X$ whose fiber at the base point $B/B \in X$ is $\C_{-\la}$.  
For each $w \in \Sym_n$, we define $X(w) := \overline{B w B / B} \subset X$ and refer to it as the Schubert variety in $X$ associated with $w$.  
We denote by $\cO_{X(w)}(\la)$ the restriction of $\cO_X(\la)$ to $X(w)$.

Using Lemma~\ref{lem:rootid}, we define a $(G \times \Gm)$-equivariant vector subbundle
\[
T_{\Psi}^* X := G \times^B \gn(\Psi) \subset G \times^B \gn \cong T^*X,
\]
for a root ideal $\Psi \subset \Delta^+$, where $\Gm$ acts by fiberwise scalar dilation. Let $\pi_\Psi \colon T^*_\Psi X \to X$ denote the natural projection. For each $w \in \Sym_n$, we set
\[
T^*_\Psi X(w) := \pi_\Psi^{-1}(X(w)).
\]
We denote the restriction of $\pi_\Psi$ to $T^*_\Psi X(w)$ again by $\pi_\Psi$, by slight abuse of notation.

For a sequence $\bi := ( i_1, i_2, \ldots, i_{\ell} )$ of elements in $\tI_\af$, we define the associated $\tbB$-schemes by
\begin{equation}
\widetilde{Z}(\bi) := \tbP_{i_1} \times^{\tbB} \tbP_{i_2} \times^{\tbB} \cdots \times^{\tbB} \tbP_{i_{\ell}},
\qquad
Z(\bi) := \widetilde{Z}(\bi) / \tbB.
\label{eqn:defZi}
\end{equation}
By convention, we set $Z(\varnothing) := \mathrm{pt}$.

\begin{lem}[Kumar {\cite[\S7.1]{Kum02}}]\label{lem:demseq}
Let $\bi := ( i_1, i_2,\ldots, i_{\ell} )$ be a sequence of elements in $\tI_\af$. Then the following statements hold:
\begin{enumerate}
\item Let $\bi^{\flat}$ be the sequence obtained by omitting the last element $i_{\ell}$ of $\bi$. Then $Z(\bi)$ is a $\P^1$-fibration over $Z(\bi^{\flat})$, whose fiber is isomorphic to $\tbP_{i_\ell}/\tbB$.
\item Let $1 \le j_1 < j_2 < \cdots < j_m \le \ell$, and set $\bi' := (i_{j_1}, i_{j_2}, \ldots, i_{j_m})$. Then there is a closed $\tbB$-equivariant embedding $Z(\bi') \hookrightarrow Z(\bi)$ induced by the group homomorphism
\[
\prod_{t = 1}^m \tbP_{i_{j_t}} \ni (g_{j_t}) \mapsto (g_{j}) \in \prod_{j=1}^\ell \tbP_{i_j},
\]
where $g_j := 1 \in \tbB$ for all $j \notin \{j_1, \ldots, j_m\}$.
\end{enumerate}
\end{lem}

For any rational $\tbB$-module $M$ with finite-dimensional $\tT$-weight spaces, we define a vector bundle
\[
\cE _{\bi} ( M ) := \widetilde{Z} ( \bi ) \times ^{\tbB} M^{\vee} \longrightarrow Z (\bi).
\]
In the special case where $M \cong \C_{\La}$ for some $\tT$-weight $\La$, we write $\cO_{\bi} ( \La ) := \cE_{\bi} ( \C_{\La} )$. By Lemma~\ref{lem:demseq}(2), the restriction of $\cE_{\bi}(M)$ to $Z(\bi')$ is naturally identified with $\cE_{\bi'}(M)$ as a $\tbB$-equivariant vector bundle.

\begin{defn}[Demazure functors]
The (covariant) functor that assigns to a rational $\tbB$-module $M$ the dual space $\Gamma ( Z ( \bi ), \mathcal E_{\bi} ( M ) )^{\vee}$ is called the \emph{Demazure functor} associated with the sequence $\bi$, and is denoted by $\mathscr{D}_\bi$. In particular, for $i \in \tI_\af$, we write $\mathscr{D}_i := \mathscr{D}_{(i)}$.

We also define the contragredient variant by
\[
\mathscr{D}_i^\dag ( \bullet ) := \left( \mathscr{D}_i ( \bullet^\vee ) \right)^\vee.
\]
\end{defn}

If a sequence $\bi$ in $\tI_\af$ is the concatenation of two sequences $\bi_1$ and $\bi_2$, then the corresponding Demazure functors satisfy
\[
\mathscr{D}_{\bi} \cong \mathscr{D}_{\bi_1} \circ \mathscr{D}_{\bi_2}
\]
by repeated applications of Lemma~\ref{lem:demseq}(1).

\begin{defn}
Let $\tL$ be a free abelian monoid, and let $R$ be a $\tL$-graded $\C$-algebra. We say that $R$ is \emph{$\tbB$-equivariant} if the following conditions hold:
\begin{itemize}
  \item For each $a \in \tL$, the graded component $R_a$ carries a rational $\tbB$-module structure.
  \item The multiplication maps $R_a \otimes R_b \to R_{a+b}$ are $\tbB$-equivariant.
  \item $R_0 = \C$ is equipped with the trivial $\tbB$-action.
\end{itemize}
\end{defn}

\begin{lem}\label{lem:ring-push}
Let $\tL$ be a free abelian monoid, and let $R$ be a $\tbB$-equivariant $\tL$-graded $\C$-algebra. Then for each $i \in \tI_\af$, the module $\mathscr{D}_i^\dag ( R )$ naturally inherits the structure of a $\tbB$-equivariant $\tL$-graded $\C$-algebra. Moreover, the following properties hold:
\begin{itemize}
\item If $R$ is commutative, then so is $\mathscr{D}_i^\dag ( R )$.
\item If $R$ is integral, then so is $\mathscr{D}_i^\dag ( R )$.
\item If $R$ is integrally closed, then so is $\mathscr{D}_i^\dag ( R )$.
\end{itemize}
\end{lem}

\begin{proof}
The $\tbB$-equivariant algebra $R$ determines a $\tL$-graded $\tbP_i$-equivariant sheaf of algebras $\cE_i(R^{\vee})$ over $Z(i) = \widetilde{Z}(i) / \tbB = \tbP_i/\tbB$. Consequently, its global sections form a $\tL$-graded algebra with a rational $\tbP_i$-action on each graded component, compatible with multiplication. The degree-zero part of $\mathscr{D}_i^\dag ( R )$ is given by
\[
\C = \Gamma ( Z(i), \cO_{Z(i)} ) = \Gamma ( \P^1, \cO_{\P^1} ).
\]

If $R$ is commutative, then $\cE_i(R^{\vee})$ is a sheaf of commutative algebras, and hence $\mathscr{D}_i^\dag ( R )$ is also commutative.

Now suppose $R$ is integral. Then $R \otimes_{\C} \C[t]$ remains integral. Moreover, if $R$ is integrally closed, then so is $R \otimes_{\C} \C[t]$. This can be verified inductively by examining the coefficients of $t$ in an integral dependence relation, starting from the lowest-degree term.

For each $x \in \P^1$, there exists an affine open neighborhood $U_x$ with a local coordinate $t_x$ such that
\[
\Gamma(U_x, \cE_i(R^{\vee})) \cong R \otimes_{R_0} \C[t_x].
\]
Since $\P^1 = \bigcup_{x \in \P^1} U_x$, we obtain
\[
\mathscr{D}_i^\dag ( R ) = \Gamma ( \P^1, \cE_i ( R^\vee ) ) = \bigcap_{x \in \P^1} R \otimes_{\C} \C[t_x]
\]
inside the algebra of rational sections.

It follows that $\mathscr{D}_i^\dag ( R )$ is integral if $R$ is, and integrally closed if $R$ is integrally closed. The latter follows because the intersection of integrally closed domains with a common field of fractions is again integrally closed.

This completes the proof.
\end{proof}

\begin{thm}[Joseph {\cite{Jos85}}]\label{thm:Jos}
For each $i \in \tI_\af$, the following hold:
\begin{enumerate}
\item There exists a natural transformation $\mathrm{Id} \rightarrow \mathscr{D}_i$.
\item There is an isomorphism of functors $\mathscr{D}_i \stackrel{\cong}{\longrightarrow} \mathscr{D}_i \circ \mathscr{D}_i$.
\item For any rational $\tbP_i$-module $M$, there is an isomorphism of functors
$$\mathscr{D}_i ( M \otimes \bullet ) \cong M \otimes \mathscr{D}_i ( \bullet ).$$
\item Let $w \in \tSym_n$ admit two reduced expressions $\bi$ and $\bi'$ connected by a sequence of braid relations. Then there is an isomorphism of functors $\mathscr{D}_\bi \cong \mathscr{D}_{\bi'}$.
\end{enumerate}
Moreover, the functor $\mathscr{D}_i$ maps finite-dimensional rational $\tbB$-modules to finite-dimensional rational $\tbP_i$-modules, which may be regarded as $\tbB$-modules via restriction.
\end{thm}

For $m\ge0$, we use the convention
\[
\mathbb L^{-m}\mathscr D_{\bi}(M)
:=
H^m\bigl(Z(\bi),\cE_{\bi}(M)\bigr)^\vee
\]
since $\mathscr D_{\bi}$ is right exact.

\begin{cor}\label{cor:DLcomm}
For distinct $i, j \in \tI_\af$, there are canonical isomorphisms
\begin{equation}
\mathbb L^{\bullet} \mathscr{D}_i \left( \C_{\La_j} \otimes \bullet \right) \cong \C_{\La_j} \otimes \mathbb L^{\bullet} \mathscr{D}_i ( \bullet ) \quad \text{and} \quad \mathbb L^{\bullet} \mathscr{D}_i(\C) \cong \C.\label{eqn:DLcomm}
\end{equation}
\end{cor}

\begin{proof}
If $i\not\equiv j\pmod n$, then the formula defining the affine
Weyl group action gives $s_i(\La_j)=\La_j$. Hence
$\C_{\La_j}$ extends to a rational $\tbP_i$-module on which
$\SL(2,i)$ acts trivially. The first isomorphism follows from
Theorem~\ref{thm:Jos}(3), after passing to derived functors,
since tensoring with $\C_{\La_j}$ is exact.
The second follows from $Z(i)\cong\P^1$ and
$H^{>0}(\P^1,\cO_{\P^1})=0$.
\end{proof}

By Theorem~\ref{thm:Jos}(4), the functor $\mathscr{D}_w := \mathscr{D}_\bi$ is well-defined for each $w \in \tSym_n$. We set $\mathscr{D}_w^\dag ( \bullet ) := \left( \mathscr{D}_w( \bullet^\vee ) \right)^\vee$. By Theorem~\ref{thm:Jos}(1), there exists a natural transformation $\mathscr{D}_w \to \mathscr{D}_v$ whenever $w < v$.

\begin{thm}[Demazure character formula; see, e.g.,~\cite{Kum02}]\label{thm:DCF}
The following statements hold:
\begin{enumerate}
\item Let $\la \in \sP^+$ and $w \in \Sym_n$, and fix a reduced expression $\bi$ of $w$. Then
\begin{align*}
H^m\bigl(X(w), \cO_{X(w)}(\la)\bigr)^{\vee} & \cong H^m\bigl(Z(\bi), \cO_\bi(\la)\bigr)^\vee \\
& \cong \mathbb{L}^{-m} \mathscr{D}_\bi(\C_\la) \cong
\begin{cases}
V_w(\la) & \text{if } m = 0, \\
0 & \text{otherwise}.
\end{cases}
\end{align*}
\item Let $\La \in \sP^+_\af$, and let $\bi$ be a sequence of elements in $\tI_\af$. Then there exists $w \in \tSym_n$ such that
\[
H^m\bigl(Z(\bi), \cO_\bi(\La)\bigr)^{\vee} \cong \mathbb{L}^{-m} \mathscr{D}_\bi(\C_\La) \cong
\begin{cases}
L_w(\La) & \text{if } m = 0, \\
0 & \text{otherwise}.
\end{cases}
\]
\item The line bundle $\cO_\bi(\La)$ on $Z(\bi)$ is base-point-free for each $\La \in \sP^+_\af$.
\end{enumerate}
\end{thm}

\begin{proof}
Assertions~(1) and~(2) are special cases of~\cite[Corollary~8.1.26]{Kum02}, while~(3) follows from~(2) and~\cite[Proposition~7.1.15]{Kum02}.
\end{proof}

\subsection{Affine Demazure modules}\label{subsec:adm}

For each $\la \in \sP$ (regarded as an element of $\sP_\af$ via the standard lift) and $k \in \Z_{>0}$, there exists $w \in \tSym_n$ such that
\begin{equation}
\la + k \wp = w \La \in \sP_\af, \quad \La \in \sP_\af^+ \label{eqn:lL-rel}
\end{equation}
as ensured by~\cite[Corollary 10.1]{Kac}. We define the level-$k$ Demazure module by
\[
D^{(k)}_\la := \mathscr{D}_w ( \C_\La ) \cong L_w ( \La ) \subset L ( \La ).
\]
This is a finite-dimensional rational $\tbB$-module, independent of the choice of $w$ satisfying~\eqref{eqn:lL-rel}.

\begin{defn}
Let $k\in\Z_{>0}$. A finite-dimensional $\tbB$-module $M$ is said
to be \emph{$D^{(k)}$-filtered} if it has a finite filtration whose
successive subquotients are isomorphic to
$\C_{r\delta}\otimes D_\la^{(k)}$ for some $r\in\Z$ and
$\la\in\sP$. Such a filtration is called a
\emph{$D^{(k)}$-filtration}.
\end{defn}

\begin{thm}[Joseph~\cite{Jos06}; see also~\cite{Nao12,Kat26b}]\label{thm:D-branch}
Let $\la \in \sP$ and $k \in \Z_{>0}$. Then
\begin{enumerate}
  \item For each $\Theta \in \{ \La_i \}_{i \in \tI_\af} \cup \{\wp\}$, the module $D_\la^{(k)} \otimes \C_{\Theta}$ admits a $D^{(k+1)}$-filtration.
  \item If $M$ is $D^{(k)}$-filtered and $i \in \tI_\af$, then $\mathbb{L}^{<0} \mathscr{D}_i(M) = 0$, and $\mathscr{D}_i(M)$ is again $D^{(k)}$-filtered.
\end{enumerate}
\end{thm}

\begin{proof}
For~(1), choose $w$ and $\La$ as in~\eqref{eqn:lL-rel}, and apply~\cite[Theorem~5.22]{Jos06} at $q=1$ with the two dominant affine weights
$\Theta$ and $\La$ and the Weyl group element $w$. The extremal
weights of the successive subquotients have the form
$\mu_a+(k+1)\wp+r_a\delta$ ($\mu_a \in \sP, r_a \in \Z$). Hence, under the convention
in~\eqref{eqn:lL-rel}, these subquotients are isomorphic to
$\C_{r_a\delta}\otimes D_{\mu_a}^{(k+1)}$. This gives the asserted
$D^{(k+1)}$-filtration. See also~\cite[Remark~4.15]{Nao12} for the
case $n=2$ and~\cite[Theorem~10.5(2)]{Kat26b} for an alternative proof.

Since $\C_{r\delta}$ extends to every $\tbP_i$,
Theorem~\ref{thm:Jos}(3) shows that $\mathscr D_i$ commutes with
grading shifts. Assertion~(2) now follows by applying
Theorem~\ref{thm:DCF}(2) to the successive subquotients of a
$D^{(k)}$-filtration.
\end{proof}

\begin{cor}[Demazure module branching]\label{cor:D-branch}
Let $k \in \Z_{> 0}$, $i \in \tI_\af$, and $w \in \tSym_n$. If $M$ is a $D^{(k)}$-filtered module and $m \in \Z_{\ge 0}$, then we have $\mathbb{L}^{<0} \mathscr{D}_w ( \C_{m \La_i} \otimes M ) = 0$, and the resulting $\tbB$-module $\mathscr{D}_w ( \C_{m \La_i} \otimes M )$ is $D^{(m+k)}$-filtered. Moreover, there is a natural inclusion
\begin{equation}
\C_{m \La_i} \otimes M \subset \mathscr{D}_w ( \C_{m \La_i} \otimes M ). \label{eqn:Dtinc}
\end{equation}
\end{cor}

\begin{proof}
Let $N$ be a finite-dimensional $\tT$-semisimple $\tbB$-module fitting into a short exact sequence
\[
0 \rightarrow N_1 \rightarrow N \rightarrow N_2 \rightarrow 0,
\]
where $N_2$ is a Demazure module and $N_1 \subset \mathscr{D}_w ( N_1 )$. Applying the Leray spectral sequence for $\mathbb{L}^{\bullet}\mathscr{D}_{\bi}$, with $\bi$ a reduced expression of $w$, we obtain
\[
\mathbb{L}^{<0} \mathscr{D}_w ( N_2 ) = 0
\]
by Theorem~\ref{thm:D-branch}(2).

Consider the following commutative diagram of short exact sequences:
\begin{equation}
\xymatrix{
 & 0 \ar[r] & N_1 \ar[r] \ar@{^{(}->}[d] & N \ar[r] \ar[d] & N_2 \ar[d]^{\imath} \ar[r] & 0 \\
0 \ar@{=}[r] & \mathbb{L}^{-1} \mathscr{D}_{w}(N_2) \ar[r] & \mathscr{D}_{w}(N_1) \ar[r] & \mathscr{D}_{w}(N) \ar[r] & \mathscr{D}_{w}(N_2) \ar[r] & 0
}.
\label{eqn:comm}
\end{equation}
The map $\imath$ is injective by Theorem~\ref{thm:DCF}(2) and the inclusion relations of Demazure modules. Thus, by the five lemma, we deduce that $N \subset \mathscr{D}_w ( N )$.

Suppose further that $\mathbb{L}^{<0} \mathscr{D}_w(N_1) = 0$.  
Then the long exact sequence associated with the bottom row of~\eqref{eqn:comm} implies that $\mathbb{L}^{<0} \mathscr{D}_w(N) = 0$.

We now apply Theorem~\ref{thm:D-branch}(1) iteratively $m$ times to obtain a $D^{(m+k)}$-filtration on $\C_{m \La_i} \otimes M$. The assertion then follows by induction on the length of the filtration, using the arguments above.
\end{proof}

\section{The rotation theorem and its interpretation}\label{sec:rot}

We divide this section into two parts. The first lifts the
Blasiak--Morse--Pun rotation formula to Demazure functors and defines the
module $M_w^\Psi(\la)$ in~\eqref{eqn:defM}. The second introduces the
module $N_w^\Psi(\la)$ in~\eqref{eqn:defN}, whose presentation is
adapted to the geometric construction beginning in
Section~\ref{sec:XP}. Proposition~\ref{prop:coinc} identifies these two
modules whenever $w$ is $\Psi$-tame.

We retain the setting of Section~\ref{sec:prelim}. For $\la\in\Par$, write
\[
m_i(\la):=
\begin{cases}
\la_i-\la_{i+1} & \text{if }1\le i<n,\\
\la_n & \text{if }i=n.
\end{cases}
\]
Thus, $\la=\sum_{i=1}^n m_i(\la)\varpi_i$.

For $1\le i<n$ and $0\le e\le n+i$, set
\[
\mathsf{L}(i,e):=
\begin{cases}
 i-e & \text{if }0\le e<i,\\
 i+n-e & \text{if }i\le e\le n+i.
\end{cases}
\]
We define
\[
\mathscr C_{i,e}
:=
\mathscr D_{i-1}\circ
\mathscr D_{i-2}\circ\cdots\circ
\mathscr D_{i-\mathsf{L}(i,e)},
\]
where the composition has exactly $\mathsf{L}(i,e)$ factors and all indices
of the Demazure functors are read modulo $n$. If $\mathsf{L}(i,e)=0$, the
composition is the identity functor. For $0\le e<i$, we have
\[
\mathsf L(i,e)=\mathsf L(i,e+n),
\quad\text{and hence} \quad
\mathscr C_{i,e}=\mathscr C_{i,e+n}.
\]
We use this convention when writing strings of second indices
cyclically.

For $1\le i<n$ and $1\le e\le n$, we also incorporate the
character twist corresponding to the right endpoint by setting
\begin{equation}
\mathscr C_{i,e}(\la)(\bullet)
:=
\mathscr C_{i,e}
\bigl(\C_{m_e(\la)\La_e}\otimes\bullet\bigr).
\label{eqn:defCie}
\end{equation}

\subsection*{The functorial form of the rotation theorem}

The following module is the lift of the Blasiak--Morse--Pun operator formula
to Demazure functors. For a root ideal $\Psi\subset\Delta^+$, $w\in\Sym_n$,
and $\la\in\Par$, define
\begin{align}
M^\Psi_w(\la)
:={}&\mathscr D_w\Bigl(
\C_{m_1(\la)\La_1}\otimes
\mathscr C_{1,e_1(\Psi)}\bigl(
\C_{m_2(\la)\La_2}\otimes
\mathscr C_{2,e_2(\Psi)}\bigl(
\cdots \nonumber\\
&\hspace{20mm}
\C_{m_{n-1}(\la)\La_{n-1}}\otimes
\mathscr C_{n-1,e_{n-1}(\Psi)}
\bigl(\C_{m_n(\la)\La_n}\bigr)
\cdots
\bigr)
\bigr)
\Bigr).
\label{eqn:defM}
\end{align}

We first record that the iterated functors in~\eqref{eqn:defM} are acyclic.
This allows us to pass from the functorial expression to the corresponding
operator formula for graded characters without correction terms.

\begin{prop}\label{prop:inclDM}
Let $\Psi\subset\Delta^+$ be a root ideal, and let $w\in\Sym_n$ and
$\la\in\Par$. Then the iterated expression on the right-hand side
of~\eqref{eqn:defM} is acyclic; equivalently, all its negative total derived
functors vanish.
\end{prop}

\begin{proof}
Let $i\in\tI_\af$, $m\in\Z_{\ge0}$, and $v\in\tSym_n$. For a
$D^{(k)}$-filtered $\tbB$-module $Q$, the relevant Leray spectral sequence is
\[
\mathbb L^r\mathscr D_i
\bigl(\C_{m\La_i}\otimes\mathbb L^s\mathscr D_v(Q)\bigr)
\Rightarrow
\mathbb L^{r+s}
\bigl(\mathscr D_i\circ
(\C_{m\La_i}\otimes\mathscr D_v)\bigr)(Q).
\]
By Corollary~\ref{cor:D-branch}, the inner module is acyclic and its
degree-zero image is again Demazure-filtered. Applying the same corollary to
the outer functor shows that the $E_2$-page is concentrated in bidegree
$(0,0)$. Iterating this argument from the rightmost factor
of~\eqref{eqn:defM} proves the assertion.
\end{proof}

\begin{thm}[Rotation theorem; Blasiak--Morse--Pun~{\cite[Theorem~2.3]{BMP}}]
\label{thm:BMP}
Let $\Psi\subset\Delta^+$ be a root ideal, and suppose that
$w\in\Sym_n$ is $\Psi$-tame. Then, for every $\la\in\Par$, we have
\[
H(\Psi;\la;w)
=
\left[\gch M^\Psi_w(\la)\right]_{q\mapsto q^{-1}},
\]
where $H(\Psi;\la;w)$ is the tame nonsymmetric Catalan function defined
in~\cite[(2.2)]{BMP}.
\end{thm}

\begin{rem}\label{rem:rotation-conventions}
To relate the statement above to the operator formula in~\cite{BMP}, let
$\Phi$ be the lift of the cyclic automorphism of the affine Dynkin diagram
used in~\cite[(2.4)]{BMP}. It satisfies
\[
\Phi\circ\pi_i=\pi_{i+1}\circ\Phi
\qquad(0\le i<n),
\]
where $\pi_i$ is the Demazure operator corresponding to $\mathscr D_i$ and
$\pi_0:=\pi_n$. Moving every occurrence of $\Phi$ to the right in
\cite[(2.5)]{BMP} yields the graded character of~\eqref{eqn:defM};
Proposition~\ref{prop:inclDM} ensures that no higher derived terms occur.
Moreover, $\Phi(x_n)=qx_1$ in the convention of~\cite{BMP}, whereas our convention gives $\Phi(X_n)=q^{-1}X_1$. This accounts for the
substitution $q\mapsto q^{-1}$ in Theorem~\ref{thm:BMP}.

Since $w_0$ is $\Psi$-tame for every root ideal $\Psi$, the theorem includes
the ordinary Catalan functions $H(\Psi;\la)=H(\Psi;\la;w_0)$. In particular,
it includes their $k$-Schur specializations.
\end{rem}

\subsection*{The geometric presentation}

We next reorganize the functors in~\eqref{eqn:defM} according to the column
blocks of $\Psi$. For $1\le j\le\ell(\Psi)$, define
\[
\mathscr C_j^\Psi(\la)
:=
\mathscr C_{\ti_j(\Psi),\te_j(\Psi)}(\la)
\circ
\mathscr C_{\ti_j(\Psi),\te_j(\Psi)+1}(\la)
\circ\cdots\circ
\mathscr C_{\ti_j(\Psi),\te_{j+1}(\Psi)-1}(\la).
\]
We also set
\[
\la(\Psi):=
\sum_{a=1}^{d_1(\Psi)}m_a(\la)\La_a
\]
and define
\begin{equation}
N^\Psi_w(\la)
:=
\mathscr D_w\left(
\C_{\la(\Psi)}\otimes
\bigl(
\mathscr C_1^\Psi(\la)\circ
\mathscr C_2^\Psi(\la)\circ\cdots\circ
\mathscr C_{\ell(\Psi)}^\Psi(\la)
\bigr)(\C)
\right).
\label{eqn:defN}
\end{equation}

The ordering in~\eqref{eqn:defN} is chosen to optimize the functorial
presentation. It is precisely this optimized form, rather than the
ordering in~\eqref{eqn:defM}, that makes possible the geometric
reinterpretation developed in the next section. For example, if
$e_{n-1}(\Psi)=n+1$, then the innermost functor
$\mathscr C_{n-1,e_{n-1}(\Psi)}$ in~\eqref{eqn:defM}
acts trivially on $\C_{m_n(\la)\La_n}$ and hence does not occur
in~\eqref{eqn:defN}.

Choose a reduced expression for $w$. For every
$1\le j\le\ell(\Psi)$ and
$\te_j(\Psi)\le e<\te_{j+1}(\Psi)$, take the consecutive
reduced word underlying $\mathscr C_{\ti_j(\Psi),e}$.
Concatenating these words in the order prescribed by~\eqref{eqn:defN},
together with the chosen reduced expression for $w$, gives a sequence
\begin{equation}
\bi=\bi(\Psi,w)=(i_1,i_2,\ldots,i_r)\label{eqn:defbiPw}
\end{equation}
in $\tI_{\af}$. This sequence is independent of $\la$.

To define the corresponding line bundle, regard
$\C_{\la(\Psi)}$ in~\eqref{eqn:defN} as the tensor product of
the characters $\C_{m_a(\la)\La_a}$ for $1\le a\le d_1(\Psi)$.
Thus, for every $1\le a\le n$, the expression~\eqref{eqn:defN}
contains a distinguished character factor
$\C_{m_a(\la)\La_a}$. Let $r_a$ be the number of Demazure functors
appearing to the left of this factor in~\eqref{eqn:defN}, and set
\[
\bi[a]:=(i_1,i_2,\ldots,i_{r_a}).
\]
Iterating Lemma~\ref{lem:demseq}(1), we obtain the forgetful morphism
associated with this initial subsequence:
\[
p_a\colon Z(\bi)\longrightarrow Z(\bi[a]).
\]
We define
\begin{equation}\label{eqn:defL-Psi}
\mathcal L^\Psi_w(\la)
:=
\bigotimes_{a=1}^n
p_a^*\cO_{\bi[a]}(\La_a)^{\otimes m_a(\la)}.
\end{equation}
By construction, the successive associated-bundle interpretation of
Demazure functors identifies the total derived functor
in~\eqref{eqn:defN} with the dual cohomology of
$\mathcal L^\Psi_w(\la)$. In particular,
Proposition~\ref{prop:inclD} below gives
\begin{equation}
N^\Psi_w(\la)^\vee
\cong H^0\bigl(Z(\bi),\mathcal L^\Psi_w(\la)\bigr),
\qquad
H^{>0}\bigl(Z(\bi),\mathcal L^\Psi_w(\la)\bigr)=0.
\label{eqn:N-BS-sections}
\end{equation}
Each factor in~\eqref{eqn:defL-Psi} is base-point-free by
Theorem~\ref{thm:DCF}(3). Hence
$\mathcal L^\Psi_w(\la)$ is base-point-free and, in particular,
$N^\Psi_w(\la)\neq0$.

For later use, we also introduce the partial geometric presentations.
Let $e_1(\Psi)\le k\le n$, and choose the unique
$1\le j\le\ell(\Psi)$ such that
$\te_j(\Psi)\le k<\te_{j+1}(\Psi)$ (see Lemma~\ref{lem:basicinv}). We set
\begin{align}\nonumber
N^\Psi(\la;k)
:={}&
\Bigl(
\mathscr C_{\ti_j(\Psi),k}(\la)\circ\cdots\circ
\mathscr C_{\ti_j(\Psi),\te_{j+1}(\Psi)-1}(\la) \circ
\mathscr C_{j+1}^\Psi(\la)\circ\cdots\circ
\mathscr C_{\ell(\Psi)}^\Psi(\la)
\Bigr)(\C)\\\label{eqn:N(k)}
={}&
\bigl(
\mathscr C_{\hs_k(\Psi),k}(\la)\circ
\mathscr C_{\hs_{k+1}(\Psi),k+1}(\la) \circ\cdots\circ
\mathscr C_{\hs_n(\Psi),n}(\la)
\bigr)(\C),
\end{align}
where the second equality follows from~\eqref{eqn:hk-piecewise}.

\begin{prop}\label{prop:inclD}
Let $\Psi\subset\Delta^+$ be a root ideal, and let $w\in\Sym_n$ and
$\la\in\Par$. Then the total complex associated with~\eqref{eqn:defN}
satisfies
\[
\mathbb L^{<0}\left(
\mathscr D_w\left(
\C_{\la(\Psi)}\otimes
\bigl(
\mathscr C_1^\Psi(\la)\circ\cdots\circ
\mathscr C_{\ell(\Psi)}^\Psi(\la)
\bigr)(\C)
\right)
\right)=0.
\]
Moreover, for every $e_1(\Psi)\le k\le n$, we have
\[
\mathbb L^{<0}\left(
\bigl(
\mathscr C_{\hs_k(\Psi),k}(\la)\circ
\mathscr C_{\hs_{k+1}(\Psi),k+1}(\la)\circ\cdots\circ
\mathscr C_{\hs_n(\Psi),n}(\la)
\bigr)(\C)
\right)=0.
\]
\end{prop}

\begin{proof}
Both assertions follow by iterating the Leray spectral sequence and applying
Corollary~\ref{cor:D-branch}, exactly as in the proof of
Proposition~\ref{prop:inclDM}.
\end{proof}

The arguments in Propositions~\ref{prop:inclDM} and~\ref{prop:inclD} are also closely related to the Bott--Samelson vanishing argument of~\cite{LLM02}; the authors note that their results extend to symmetrizable Kac--Moody groups.

The following comparison is the main result of this section. It converts the
Blasiak--Morse--Pun presentation into the geometric presentation used in
the rest of the paper.

\begin{prop}\label{prop:coinc}
Let $\Psi\subset\Delta^+$ be a root ideal, and suppose that
$w\in\Sym_n$ is $\Psi$-tame. Then, for every $\la\in\Par$, there is an
isomorphism of $\tbB$-modules
\[
N^\Psi_w(\la)\cong M^\Psi_w(\la).
\]
\end{prop}

We prepare the proof by recording two consequences of the braid relations for
Demazure functors.

\begin{lem}\label{lem:ci-move}
Let $1\le i<e\le n$, and suppose that
$M\cong\mathscr D_j(M)$ for all $0<j<i$. Then
\[
\mathscr C_{i,e}(M) \cong
\mathscr D_j\bigl(\mathscr C_{i,e}(M)\bigr)
\qquad(0\le j<i).
\]
\end{lem}

\begin{proof}
Let $v' \in \Sym_i$ and $v \in \Sym_{i+1}$ denote the longest elements in the subgroups $\langle s_{i-1}, \ldots, s_1 \rangle$ and $\langle s_{i-1}, \ldots, s_0 \rangle \subset \tSym_n$, respectively. By assumption, we have $M \cong \mathscr{D}_{v'}(M)$, and hence $\mathscr{C}_{i,e}(M) \cong \mathscr{C}_{i,e}(\mathscr{D}_{v'}(M))$.

Since $v s_j < v$ for all $0 \le j < i$, we compute:
\[
\mathscr{C}_{i,e} \circ \mathscr{D}_{v'} \cong \left( \mathscr{D}_{i-1} \circ \cdots \circ \mathscr{D}_0 \right) \circ \left( \mathscr{D}_{n-1} \circ \cdots \circ \mathscr{D}_e \right) \circ \mathscr{D}_{v'} \cong \mathscr{D}_v \circ \mathscr{C}_{i,e},
\]
where the second isomorphism follows from Theorem~\ref{thm:Jos}(2)(4). Since $s_j v < v$ for all $0 \le j < i$, it follows that $\mathscr{D}_v \cong \mathscr{D}_j \circ \mathscr{D}_v$ by the same theorem. Therefore, $\mathscr{C}_{i,e}(M) \cong \mathscr{D}_j(\mathscr{C}_{i,e}(M))$.
\end{proof}

\begin{cor}\label{cor:ci-move}
Let $1\le i<e\le e'\le n$. If
$M\cong\mathscr D_j(M)$ for all $0\le j<i$ and $e'<j<n$, then
\[
\mathscr C_{i,e}(M)
\cong
\mathscr D_j\bigl(\mathscr C_{i,e}(M)\bigr)
\qquad
(0\le j<i\ \text{or}\ e'\le j<n).
\]
\end{cor}

\begin{proof}
The affine Dynkin diagram of type $\mathsf{A}_{n-1}^{(1)}$ admits an automorphism given by cyclically rotating the indices of the simple roots.  
Applying this automorphism, specifically by adding $(n - e')$ modulo $n$ to all indices, reduces the claim to the case treated in Lemma~\ref{lem:ci-move}.
\end{proof}

\begin{lem}\label{lem:b-move}
Let $1\le i<e\le n$. For $e\le j<n$ or $0\le j<i-1$, we have
\[
\mathscr D_j\circ\mathscr C_{i,e}
\cong
\mathscr C_{i,e}\circ\mathscr D_{j+1}.
\]
\end{lem}

\begin{proof}
The claim follows from the braid relation
\[
\mathscr{D}_j \circ \mathscr{D}_{j+1} \circ \mathscr{D}_j 
\cong 
\mathscr{D}_{j+1} \circ \mathscr{D}_j \circ \mathscr{D}_{j+1},
\]
which is a special case of Theorem~\ref{thm:Jos}(4).  
It suffices to observe that the remaining functors in $\mathscr{C}_{i,e}$ commute with $\mathscr{D}_j$ and $\mathscr{D}_{j+1}$, respectively:  
namely, $\mathscr{D}_j$ commutes with $\mathscr{D}_{i-1}, \ldots, \mathscr{D}_{j+2}$,  
and $\mathscr{D}_{j+1}$ commutes with $\mathscr{D}_{j-1}, \ldots, \mathscr{D}_e$.
\end{proof}

\begin{cor}\label{cor:b-move}
Let $2\le i<e<n$. For $e\le e'<n$ or $0\le e'<i-1$, we have
\[
\mathscr C_{i-1,e'}\circ\mathscr C_{i,e}
\cong
\mathscr C_{i,e}\circ\mathscr C_{i,e'+1}.
\]
\end{cor}

\begin{proof}
We apply Lemma~\ref{lem:b-move} to the composition  
$\mathscr{C}_{i-1,e'} \circ \mathscr{C}_{i,e} = \mathscr{D}_{i-2} \circ \cdots \circ \mathscr{D}_{e'} \circ \mathscr{C}_{i,e}$,  
and move each $\mathscr{D}_j$ past $\mathscr{C}_{i,e}$ using Lemma~\ref{lem:b-move} repeatedly.  
This yields
\[
\mathscr{D}_{i-2} \circ \cdots \circ \mathscr{D}_{e'} \circ \mathscr{C}_{i,e} 
\cong 
\mathscr{C}_{i,e} \circ \mathscr{D}_{i-1} \circ \cdots \circ \mathscr{D}_{e'+1},
\]
from which the desired isomorphism follows.
\end{proof}

Before proving Proposition~\ref{prop:coinc}, we illustrate the block
transformation in the example introduced in Section~\ref{sec:prelim}.

\begin{ex}\label{ex:n6-re}
We employ the setting of Example~\ref{ex:n6}, taking $\la = \varpi_n$.  

The presentation of $N_w^\Psi(\varpi_n)$ contains the composition
\[
\mathscr{C}_{2,3} \circ \mathscr{C}_{2,4} \circ \mathscr{C}_{2,5} \circ \mathscr{C}_{3,6} 
= (\mathscr{D}_1 \mathscr{D}_0 \mathscr{D}_5 \mathscr{D}_4 \mathscr{D}_3)
  (\mathscr{D}_1 \mathscr{D}_0 \mathscr{D}_5 \mathscr{D}_4)
  (\mathscr{D}_1 \mathscr{D}_0 \mathscr{D}_5)
  (\mathscr{D}_2 \mathscr{D}_1 \mathscr{D}_0).
\]
The corresponding part of the presentation of $M_w^\Psi(\varpi_n)$ is
\begin{equation}
(\mathscr{D}_0 \mathscr{D}_5 \mathscr{D}_4 \mathscr{D}_3)
(\mathscr{D}_1 \mathscr{D}_0 \mathscr{D}_5 \mathscr{D}_4 \mathscr{D}_3)
(\mathscr{D}_2 \mathscr{D}_1 \mathscr{D}_0)
(\mathscr{D}_3 \mathscr{D}_2 \mathscr{D}_1)
(\mathscr{D}_4 \mathscr{D}_3 \mathscr{D}_2 \mathscr{D}_1).
\label{eqn:DC}
\end{equation}
In the full presentation of $M_w^\Psi(\varpi_n)$,
we apply~\eqref{eqn:DC} to $\C_{\La_n}$ and then let
$\mathscr D_3$, $\mathscr D_4$, and $\mathscr D_5$
act on the result from the left.

Now we transform these two presentations. Note that~\eqref{eqn:DC} simplifies to
\begin{equation}
(\mathscr{D}_0 \mathscr{D}_5 \mathscr{D}_4 \mathscr{D}_3)(\mathscr{D}_1 \mathscr{D}_0 \mathscr{D}_5 \mathscr{D}_4 \mathscr{D}_3)(\mathscr{D}_2 \mathscr{D}_1 \mathscr{D}_0),\label{eqn:DC2}
\end{equation}
since $\mathscr{D}_i(\C_{\La_n}) = \C_{\La_n}$ for all $i \ne 0$.

Here, we use the isomorphism of functors
\[
\mathscr{D}_i (\mathscr{D}_1 \mathscr{D}_0 \mathscr{D}_5 \mathscr{D}_4 \mathscr{D}_3) 
= (\mathscr{D}_1 \mathscr{D}_0 \mathscr{D}_5 \mathscr{D}_4 \mathscr{D}_3)\mathscr{D}_{i+1}
\quad \text{for } i = 3, 4, 5, 0,
\]
to transform~\eqref{eqn:DC2} into
\begin{equation}
(\mathscr{D}_1 \mathscr{D}_0 \mathscr{D}_5 \mathscr{D}_4 \mathscr{D}_3)(\mathscr{D}_1 \mathscr{D}_0 \mathscr{D}_5 \mathscr{D}_4)(\mathscr{D}_2 \mathscr{D}_1 \mathscr{D}_0).\label{eqn:DC3}
\end{equation}

Furthermore, for $i = 3, 4, 5$, we have
\[
\mathscr{D}_i (\mathscr{D}_1 \mathscr{D}_0 \mathscr{D}_5 \mathscr{D}_4 \mathscr{D}_3)(\mathscr{D}_1 \mathscr{D}_0 \mathscr{D}_5 \mathscr{D}_4) 
= (\mathscr{D}_1 \mathscr{D}_0 \mathscr{D}_5 \mathscr{D}_4 \mathscr{D}_3)(\mathscr{D}_1 \mathscr{D}_0 \mathscr{D}_5 \mathscr{D}_4)\mathscr{D}_{i+2}.
\]

Applying the left actions of $\mathscr{D}_3$, $\mathscr{D}_4$, and $\mathscr{D}_5$ to~\eqref{eqn:DC3}, we recover the identity
\[
(\mathscr{D}_1 \mathscr{D}_0 \mathscr{D}_5 \mathscr{D}_4 \mathscr{D}_3)
(\mathscr{D}_1 \mathscr{D}_0 \mathscr{D}_5 \mathscr{D}_4)
(\mathscr{D}_1 \mathscr{D}_0 \mathscr{D}_5)
(\mathscr{D}_2 \mathscr{D}_1 \mathscr{D}_0)
= \mathscr{C}_{2,3} \circ \mathscr{C}_{2,4} \circ \mathscr{C}_{2,5} \circ \mathscr{C}_{3,6}.
\]
Thus, we have $M_w^\Psi(\varpi_n) = N_w^\Psi(\varpi_n)$ in this case.
\end{ex}

\begin{proof}[Proof of Proposition~\ref{prop:coinc}]
If $\Psi=\varnothing$, then $e_a(\Psi)=n+1$ for every $a<n$
and $d_1(\Psi)=n$. Hence all the functors
$\mathscr C_{a,n+1}$ in~\eqref{eqn:defM} act trivially on the
character twists to their right, and therefore, for every $w\in\Sym_n$,
\[
M_w^\Psi(\la)
\cong
\mathscr D_w\left(
\C_{\sum_{a=1}^n m_a(\la)\La_a}
\right)
=
N_w^\Psi(\la).
\]

Henceforth, assume that $\Psi\ne\varnothing$.
Using Corollary~\ref{cor:DLcomm}, we first move every character twist in
\eqref{eqn:defM} to the left across the Demazure functors that fix it. After
discarding these trivial functors, we obtain
\begin{equation}
M^\Psi_w(\la)
\cong
\mathscr D_w\left(
\C_{\mu_0}\otimes
\mathscr C_{1,e_1(\Psi)}\left(
\C_{\mu_1}\otimes
\mathscr C_{2,e_2(\Psi)}\left(
\cdots
\mathscr C_{r,e_r(\Psi)}(\C_{\mu_r})
\cdots
\right)
\right)
\right),
\label{eqn:moddefM}
\end{equation}
where Lemma~\ref{lem:basicinv} gives
\[
r =\max\{1\le a<n\mid e_a(\Psi) = a+d_a(\Psi)\le n\}
  =\ti_{\ell(\Psi)}(\Psi).
\]
For
$0\le a\le r$, set
\[
\mu_a:=
\sum_{i=e_a(\Psi)}^{e_{a+1}(\Psi)-1}m_i(\la)\La_i.
\]
Note that $e_0(\Psi)=1$ and $e_{r+1}(\Psi)=n+1$. For the remainder of the proof, we suppress $\Psi$ from
$\ti_j(\Psi)$ and $\te_j(\Psi)$ and read all strings of affine
indices cyclically modulo $n$.

\smallskip
\noindent\emph{The block transformation.}
For $1\le j\le\ell(\Psi)$, using~\eqref{eqn:teblock}, we isolate
the following block in~\eqref{eqn:moddefM}:
\begin{equation}
\C_{\mu_{\ti_{j-1}}}\otimes
\bigl(
\mathscr C_{\ti_{j-1}+1,\te_j}\circ
\mathscr C_{\ti_{j-1}+2,\te_j}\circ\cdots\circ
\mathscr C_{\ti_j,\te_j}
\bigr)
(\C_{\mu_{\ti_j}}\otimes\bullet)
\label{eqn:parts}
\end{equation}
We claim that it can be replaced by
\begin{equation}
\C_{\mu_{\ti_{j-1}}}\otimes
\bigl(
\mathscr C_{\ti_j,\te_j}\circ
\mathscr C_{\ti_j,\te_j+1}\circ\cdots\circ
\mathscr C_{\ti_j,\ti_j-1}
\bigr)
(\C_{\mu_{\ti_j}}\otimes\bullet).
\label{eqn:mparts}
\end{equation}

For $j=1$, the $\Psi$-tameness of $w$ permits us to append
\begin{equation}
\mathscr D_{\ti_{j-1}-1},
\mathscr D_{\ti_{j-1}-2},\ldots,
\mathscr D_{\te_j}
\label{eqn:DinvL}
\end{equation}
to $\mathscr D_w$ without changing~\eqref{eqn:moddefM}. For $j>1$, the same
invariance is supplied by the preceding transformed block, as explained
below. By Lemma~\ref{lem:basicinv}, every $\La_a$ occurring in
$\mu_{\ti_{j-1}}$ satisfies
\[
\te_{j-1}\le a<\te_j.
\]
The indices in~\eqref{eqn:DinvL} lie outside this interval. Hence,
Corollary~\ref{cor:DLcomm} allows us to move these functors across
$\C_{\mu_{\ti_{j-1}}}$ and insert each of~\eqref{eqn:DinvL} immediately before
$\mathscr C_{\ti_{j-1}+1,\te_j}$.

Repeated applications of Lemma~\ref{lem:b-move} move the inserted functors
through the block and replace them by
\begin{equation}
\mathscr D_{\ti_j-1},
\mathscr D_{\ti_j-2},\ldots,
\mathscr D_{\te_j+\ti_j-\ti_{j-1}}
\label{eqn:DinvL2}
\end{equation}
immediately after $\mathscr C_{\ti_j,\te_j}$. Consequently, we may insert
\begin{equation}
\mathscr C_{\ti_j,\te_j+\ti_j-\ti_{j-1}},
\ldots,
\mathscr C_{\ti_j,\ti_j-1}
\label{eqn:DinvL3}
\end{equation}
at the same position. Corollary~\ref{cor:b-move}, applied repeatedly, gives
\[
\mathscr C_{\ti_{j-1}+1,\te_j}\circ\cdots\circ
\mathscr C_{\ti_j,\te_j}
\cong
\mathscr C_{\ti_j,\te_j}\circ
\mathscr C_{\ti_j,\te_j+1}\circ\cdots\circ
\mathscr C_{\ti_j,\te_j+\ti_j-\ti_{j-1}-1}.
\]
Together with~\eqref{eqn:DinvL3}, this proves the claimed transformation
from~\eqref{eqn:parts} to~\eqref{eqn:mparts}.

The product of the functors in~\eqref{eqn:mparts} is the Demazure functor
associated with the longest element of the finite parabolic subgroup
\begin{equation}
\langle
s_{\ti_j-1},\ldots,s_0,s_{n-1},\ldots,s_{\te_j}
\rangle\subset\tSym_n.
\label{eqn:j-blockW}
\end{equation}
It is therefore invariant under the simple Demazure functors belonging to
this subgroup. When $j<\ell(\Psi)$, we may in particular append
$\mathscr D_{\ti_j-1},\ldots,\mathscr D_{\te_{j+1}}$ after the block. These
are precisely the functors required in~\eqref{eqn:DinvL} for the next value
of $j$. This completes the induction.

\smallskip
\noindent\emph{Assembly of the blocks.}
Replace each occurrence of~\eqref{eqn:parts} in~\eqref{eqn:moddefM} by
\eqref{eqn:mparts}. For $j<\ell(\Psi)$, the excess tail is
\begin{equation}
\mathscr C_{\ti_j,\te_{j+1}},
\mathscr C_{\ti_j,\te_{j+1}+1},\ldots,
\mathscr C_{\ti_j,\ti_j-1}.
\label{eqn:C-extra}
\end{equation}
It commutes with $\C_{\mu_{\ti_j}}$ and can be absorbed into the next block;
there it acts through simple reflections belonging to the parabolic subgroup
in~\eqref{eqn:j-blockW} for $j+1$ and may be discarded. For
$j=\ell(\Psi)$, the corresponding cyclic tail is
$\mathscr C_{\ti_j,1},\ldots,\mathscr C_{\ti_j,\ti_j-1}$; it acts trivially
on the remaining character twists by Corollary~\ref{cor:DLcomm} and may also
be discarded.

The $j$-th block is consequently reduced to
\[
\C_{\mu_{\ti_{j-1}}}\otimes
\bigl(
\mathscr C_{\ti_j,\te_j}\circ
\mathscr C_{\ti_j,\te_j+1}\circ\cdots\circ
\mathscr C_{\ti_j,\te_{j+1}-1}
\bigr)
(\C_{\mu_{\ti_j}}\otimes\bullet),
\]
which is exactly the factor $\mathscr C_j^\Psi(\la)$ occurring
in~\eqref{eqn:defN}. Performing this replacement for all $j$ transforms
$M_w^\Psi(\la)$ into $N_w^\Psi(\la)$ and proves the proposition.
\end{proof}

\begin{rem}\label{rem:tame-scope}
The distinction between the tame and non-tame cases is now transparent. The
modules $M_w^\Psi(\la)$ and $N_w^\Psi(\la)$ and the acyclicity statements in
Propositions~\ref{prop:inclDM} and~\ref{prop:inclD} are defined and valid for
every $w\in\Sym_n$. Tameness is used in two places: in the
Blasiak--Morse--Pun character formula, Theorem~\ref{thm:BMP}, and in the first
step of the block transformation above, where it allows the required simple
Demazure functors to be appended to $\mathscr D_w$. Without this invariance,
the present argument gives neither the identification
$M_w^\Psi(\la)\cong N_w^\Psi(\la)$ nor a character formula for
$N_w^\Psi(\la)$ in terms of $H(\Psi;\la;w)$. We make no assertion here about
such an identification in the non-tame case.
\end{rem}

We conclude with two invariance properties of the partial modules. They will
be used in the construction of $\sX_\Psi$.

\begin{lem}\label{lem:partNinv}
Let $\Psi\subset\Delta^+$ be a root ideal and $\la\in\Par$. For every
$e_1(\Psi)\le k\le n$, the $\tbB$-module $N^\Psi(\la;k)$ defined in~\eqref{eqn:N(k)} is invariant under
\[
\mathscr D_k,\ldots,\mathscr D_{n-1},
\mathscr D_0,\ldots,\mathscr D_{\hs_k(\Psi)-1}.
\]
\end{lem}

\begin{proof}
We use downward induction on $k$. For $k=n$, the functors $\mathscr D_i$ with
$1\le i<n$ act trivially on $\C_{m_n(\la)\La_n}$. Lemma~\ref{lem:ci-move},
applied with $i=\hs_n(\Psi)$ and $e=n$, therefore proves the assertion.

Assume the result for $k+1$, and put
\[
Q_{k+1}:=
\bigl(
\mathscr C_{\hs_{k+1}(\Psi),k+1}(\la)\circ\cdots\circ
\mathscr C_{\hs_n(\Psi),n}(\la)
\bigr)(\C).
\]
The induction hypothesis and
$\hs_k(\Psi)\le\hs_{k+1}(\Psi)$ imply that $Q_{k+1}$ is invariant under
$\mathscr D_{k+1},\ldots,\mathscr D_{n-1}$ and
$\mathscr D_0,\ldots,\mathscr D_{\hs_k(\Psi)-1}$. Moreover,
Corollary~\ref{cor:DLcomm} allows the twist
$\C_{m_k(\la)\La_k}$ to pass through all these functors. Applying
Corollary~\ref{cor:ci-move} with
$i=\hs_k(\Psi)$, and $e=e'=k$ now shows that
\[
N^\Psi(\la;k)=
\mathscr C_{\hs_k(\Psi),k}(\la)(Q_{k+1})
\]
has the required invariance.
\end{proof}

\begin{lem}\label{lem:tame}
Let $\Psi\subset\Delta^+$ be a root ideal, and let $w\in\Sym_n$. Then, for
every $\la\in\Par$ and $e_1(\Psi)\le i<n$, we have
\[
N_w^\Psi(\la)\cong N_{ws_i}^\Psi(\la).
\]
\end{lem}

\begin{proof}
By Lemma~\ref{lem:partNinv}, the module
$N^\Psi(\la;e_1(\Psi))$ is invariant under
$\mathscr D_i$ for $e_1(\Psi)\le i<n$. In addition,
\[
\left\langle\alpha_i,\la(\Psi)\right\rangle
=
\left\langle
\alpha_i,
\sum_{a=1}^{d_1(\Psi)}m_a(\la)\La_a
\right\rangle
=0.
\]
Corollary~\ref{cor:DLcomm} therefore allows $\mathscr D_i$ to pass through
the character twist $\C_{\la(\Psi)}$. The braid relations and idempotence in
Theorem~\ref{thm:Jos}(2)(4) then identify the actions of $\mathscr D_w$ and
$\mathscr D_{ws_i}$ on
$\C_{\la(\Psi)}\otimes N^\Psi(\la;e_1(\Psi))$, proving the assertion.
\end{proof}

\section{Construction of the variety $\sX_{\Psi}$}\label{sec:XP}

We adopt the notation and assumptions from the previous section.

\begin{lem}\label{lem:coord}
Let $\Psi \subset \Delta^+$ be a root ideal, and let $w \in \Sym_n$ and $e_1(\Psi) \le k \le n$.
Then the $\Par$-graded vector spaces
\begin{equation}
\bigoplus_{\la \in \Par} N^{\Psi}_{w}(\la)^\vee 
\qquad \text{and} \qquad 
\bigoplus_{\la \in \Par} N^{\Psi}(\la; k)^\vee
\label{eqn:pcr}
\end{equation}
admit structures of commutative $\tbB$-equivariant $\Par$-graded $\C$-algebras.
Moreover, both are integral domains and integrally closed.
\end{lem}

\begin{proof}
The character twists arising in the constructions of the modules $N^{\Psi}_{w}(\la)$ and $N^{\Psi}(\la; k)$ are additive with respect to the monoid structure on $\Par$.
Consequently, the assertion follows by repeated applications of Lemma~\ref{lem:ring-push} from the initial case $\C [\Par] \cong \C[e^{-\La_i}\mid 1 \le i \le n]$.
\end{proof}

Given a root ideal $\Psi \subset \Delta^+$ and $w \in \Sym_n$, we define the $\tbB$-scheme
\[
\sX_{\Psi}(w) := \Proj_{\Par} \left( \bigoplus_{\la \in \Par} N^{\Psi}_w(\la)^\vee \right)
\]
as a multigraded $\Proj$ over $\C$, following the construction in~\eqref{eqn:mproj}.  
For $e_1(\Psi) \le k \le n$, we denote the auxiliary spaces in the construction by
\[
\sY_{\Psi}(k) := \Proj_{\Par} \left( \bigoplus_{\la \in \Par} N^{\Psi}(\la; k)^\vee \right).
\]
By Lemma~\ref{lem:coord}, both $\sX_{\Psi}(w)$ and $\sY_{\Psi}(k)$ are integral and normal schemes.

We also set $\sY_\Psi(n+1) := \mathrm{pt}$ with trivial $\tbG$-action. With this convention,
\[
\sY_\Psi(e_1(\Psi))=\sX_\Psi(w_0^\Psi)
\]
for every root ideal $\Psi$ (see Lemma~\ref{lem:idXsX}).

The construction in~\eqref{eqn:defL-Psi} is additive in $\la$:
there are canonical isomorphisms
\[
\mathcal L^\Psi_w(\la+\mu)
\cong
\mathcal L^\Psi_w(\la)\otimes\mathcal L^\Psi_w(\mu)
\qquad (\la,\mu\in\Par).
\]
Under the identifications~\eqref{eqn:N-BS-sections}, the multiplication
maps in the first algebra in~\eqref{eqn:pcr} agree with the multiplication
of sections induced by these isomorphisms. Consequently,
\begin{equation}\label{eqn:section-ring-BS}
\bigoplus_{\la\in\Par}N^\Psi_w(\la)^\vee
\cong
\bigoplus_{\la\in\Par}
H^0\bigl(Z(\bi),\mathcal L^\Psi_w(\la)\bigr)
\end{equation}
as $\tbB$-equivariant $\Par$-graded algebras. The same construction,
applied to the truncated sequence defining $N^\Psi(\la;k)$, realizes
the second algebra in~\eqref{eqn:pcr} as a multisection ring on the
corresponding Bott--Samelson variety.

Example~\ref{ex:n4} (at the end of this section) illustrates the construction of $\sX_{\Psi}(w_0)$ in the case $n = 4$.

\begin{cor}\label{cor:XmapP}
Let $\Psi \subset \Delta^+$ be a root ideal.
Then there exist natural $\tbB$-equivariant morphisms
\begin{equation}
\sX_\Psi(w_0^\Psi) \longrightarrow \prod_{k = e_1(\Psi)}^n \P\big(N^\Psi_{w_0^\Psi}(\varpi_k)\big) \hookrightarrow \prod_{k =1}^n \P\big(L(\La_k)\big),
\label{eqn:Xprojemb}
\end{equation}
where the second map is a closed immersion.
\end{cor}

\begin{proof}
By~\eqref{eqn:N-BS-sections} and the base-point-freeness of
$\mathcal L^\Psi_{w_0^\Psi}(\varpi_k)$, the degree-$\varpi_k$
component of the multisection ring~\eqref{eqn:section-ring-BS}
defines a $\tbB$-equivariant morphism
\begin{equation}
f_k \colon
\sX_\Psi(w_0^\Psi)
\longrightarrow
\P\bigl(N^\Psi_{w_0^\Psi}(\varpi_k)\bigr)
\qquad
(e_1(\Psi)\le k\le n).\label{eqn:deffk}
\end{equation}
Taking the product of these morphisms gives the first map
in~\eqref{eqn:Xprojemb}.

By Theorem~\ref{thm:DCF}(2), we have
\begin{align*}
N^\Psi_{w_0^\Psi}(\varpi_k) & \subset L(\La_k)
\qquad (1\le k\le n),\\
N^\Psi_{w_0^\Psi}(\varpi_i) & =\C\bv_{\La_i}
\qquad (1\le i\le d_1(\Psi)).
\end{align*}
Thus the second map in~\eqref{eqn:Xprojemb} is obtained from the
linear closed immersions
\[
\P\bigl(N^\Psi_{w_0^\Psi}(\varpi_k)\bigr)
\hookrightarrow \P\bigl(L(\La_k)\bigr)
\qquad (e_1(\Psi)\le k\le n)
\]
by adjoining the fixed points $[\bv_{\La_k}]$ for
$1\le k\le d_1(\Psi)$. Hence it is a closed immersion.
\end{proof}

\begin{lem}\label{lem:idXsX}
Let $\Psi \subset \Delta^+$ be a root ideal, and let $w \in \Sym_n$ and $e_1(\Psi) \le k \le n$.
Then we have closed immersions of $\tbB$-schemes 
\[
\sY_\Psi(k) \subset \sX_\Psi(w) \subset \sX_\Psi(w_0).
\]
In particular, $\sY_\Psi(e_1(\Psi)) = \sX_\Psi(w_0^\Psi)$.
\end{lem}

\begin{proof}
By Lemma~\ref{lem:tame}, we have an identification of homogeneous coordinate rings of
$\sY_\Psi(e_1(\Psi)) = \sX_\Psi(w_0^\Psi)$,
since $N^\Psi_e(\la) = \C_{\la(\Psi)} \otimes N^\Psi(\la; e_1(\Psi))$ for all $\la \in \Par$.
The remaining closed immersions follow from surjective maps between the corresponding homogeneous coordinate rings, obtained via repeated applications of Corollary~\ref{cor:D-branch}.
\end{proof}

The graded components of the ring~\eqref{eqn:pcr} define $\tbB$-equivariant line bundles 
$\cO_{\sX_\Psi(w)}(\la)$ on $\sX_\Psi(w)$ and $\cO_{\sY_\Psi(k)}(\la)$ on $\sY_\Psi(k)$ for each $\la \in \Par$.
These line bundles extend to all $\la \in \sP$ via duality and tensor product operations.

We now define two subgroups of $\widetilde{G}(\!(z)\!)$ as follows:
\begin{align*}
\tbP(k) &:= \left< \tbP_i \mid k \le i < n \text{ or } 0 \le i < \hs_k(\Psi) \right>, \\
G(k) &:= \left< \SL(2, i) \mid k \le i < n \text{ or } 0 \le i < \hs_k(\Psi) \right>,
\end{align*}
where $e_1(\Psi) \le k \le n$.
By convention, we set $\tbP(n+1) := \tbG$ and $G(n+1) := G$.

\begin{lem}
Let $\Psi \subset \Delta^+$ be a root ideal. For each $e_1(\Psi) \le k \le n$, we have
$G(k) \cong \SL(\hs_k(\Psi) + n - k + 1)$, $\tbP(k) = G(k) \cdot \tbB$, and the group $\tbP(k)$ is proalgebraic.
In addition, there exists a split surjective homomorphism
\[
\tbP(k) \longrightarrow G(k).
\]
\end{lem}

\begin{proof}
We invoke the Dynkin diagram automorphism of type $\mathsf{A}_{n-1}^{(1)}$,
which permutes the subgroups $\SL(2, i)$ for $i \in \tI_\af$.
Applying the cyclic shift by $(n-k+1)$ modulo $n$ to the simple roots $\{\pm\alpha_i\mid k\le i<n
\text{ or }0\le i<\hs_k(\Psi)\}$,
we observe that the corresponding one-parameter subgroups generate
$\SL(\hs_k(\Psi) + n - k + 1)$ inside $G \subset \tbG$.

Each rotated subgroup $\tbP_i$ (for $k \le i < n$ or $0 \le i < \hs_k(\Psi)$) lies in $\tbG$ and defines a closed proalgebraic subgroup.
Moreover, under the rotation, its image is the standard parabolic $P_{j} \subset G$ for $1 \le j \le \hs_k(\Psi) + n - k$, and these generate $\SL(\hs_k(\Psi) + n - k + 1)$.
It follows that $\tbP(k) = G(k) \cdot \tbB$.

Finally, the projection $z \mapsto 0$ (after applying the cyclic shift) induces the desired split surjection $\tbP(k) \to G(k)$.
\end{proof}

\begin{lem}\label{lem:Ustr}
Let $\Psi \subset \Delta^+$ be a root ideal, and let $e_1(\Psi) \le k \le n$.
Then the algebraic subgroup
\[
G(k) \cap \tbP(k+1) \subset G(k)
\]
is a maximal proper parabolic subgroup, and the natural map
\begin{equation}
G(k)/\bigl( G(k) \cap \tbP(k+1) \bigr) \longrightarrow \tbP(k)/\bigl( \tbP(k) \cap \tbP(k+1) \bigr)
\label{eqn:isom}
\end{equation}
induced by the inclusion is an isomorphism of varieties.
Moreover, the variety in~\eqref{eqn:isom} is isomorphic to the projective space $\P^{\hs_k(\Psi)+n-k}$.
\end{lem}

\begin{proof}
The Iwahori subgroup $\tbB$ is stable under the Dynkin diagram automorphism of type $\mathsf{A}_{n-1}^{(1)}$,
so $G(k) \cap \tbB$ contains a Borel subgroup of $G(k)$.
It follows that $G(k) \cap \tbP(k+1)$ is a parabolic subgroup of $G(k)$.

The inclusion induces the map~\eqref{eqn:isom},
which is an isomorphism since $\tbP(k) = G(k) \cdot \tbB$.

For $k<n$, the inequality
$\hs_k(\Psi)\le\hs_{k+1}(\Psi)$ gives
\[
\SL(2,i)\subset G(k)\cap G(k+1)
\iff
k<i<n\ \text{or}\ 0\le i<\hs_k(\Psi).
\]
For $k=n$, the convention $G(n+1)=G$ instead gives
\[
\SL(2,i)\subset G(n)\cap G
\iff
1\le i<\hs_n(\Psi).
\]
In either case, $G(k)\cap\tbP(k+1)$ is a maximal proper
parabolic subgroup of $G(k)$ with semisimple Levi factor
$\SL(\hs_k(\Psi)+n-k)$.

Thus the varieties in~\eqref{eqn:isom} are both isomorphic to
$\P^{\hs_k(\Psi)+n-k}$.
\end{proof}

\begin{lem}\label{lem:ind}
Let $\Psi \subset \Delta^+$ be a root ideal, and let $e_1(\Psi) \le k \le n$.
Suppose that $M$ is a $\tbP(k+1)$-module whose restriction to $G(k+1)$ is a rational representation.
Then we have
\[
\mathscr{C}_{\hs_k(\Psi), k}(M) \cong H^0\bigl( \tbP(k) / (\tbP(k) \cap \tbP(k+1)), \mathcal F(M) \bigr)^{\vee},
\]
where $\mathcal F(M)$ denotes the vector bundle on the quotient $\tbP(k) / (\tbP(k) \cap \tbP(k+1))$ associated to the $(\tbP(k) \cap \tbP(k+1))$-module $M^{\vee}$.
\end{lem}

\begin{proof}
Let
\[
 \bi_k:=\bigl(\hs_k(\Psi)-1,\ldots,1,0,n-1,\ldots,k\bigr).
\]
This is the ordered sequence whose Demazure functor is
$\mathscr D_{\bi_k}=\mathscr C_{\hs_k(\Psi),k}$; in particular, the order is
cyclic and is not the ordinary interval from $k$ to $\hs_k(\Psi)-1$.
Multiplication of the factors in $\widetilde Z(\bi_k)$ gives a natural
projection
\[
\pi_{\bi_k} \colon Z(\bi_k) \longrightarrow \tbP(k) / (\tbP(k) \cap \tbP(k+1)) \cong \P^{\hs_k(\Psi) + n - k}.
\]
It is the Bott--Samelson resolution of the projective homogeneous space on
the right. We denote $\tbP(k) / (\tbP(k) \cap \tbP(k+1))$ by $W$.

Consider the Leray spectral sequence associated with $\pi_{\bi_k}$:
\begin{equation}
H^q\bigl( W,  \mathbb R^p (\pi_{\bi_k})_* \pi_{\bi_k}^* \mathcal F(M) \bigr) \Rightarrow H^{q+p}\bigl(Z(\bi_k), \pi_{\bi_k}^* \mathcal F(M)\bigr).\label{eqn:stepLeray}
\end{equation}

By the projection formula,
\[
\mathbb R^{\bullet} (\pi_{\bi_k})_* \pi_{\bi_k}^* \mathcal F(M) \cong \mathbb R^{\bullet} (\pi_{\bi_k})_* \cO_{Z(\bi_k)} \otimes \mathcal F(M).
\]
It is known from \cite[Theorem~8.2.2(c) and A.24]{Kum02} that $\mathbb R^{\bullet} (\pi_{\bi_k})_* \cO_{Z(\bi_k)} = \cO_{W}$, so the spectral sequence~\eqref{eqn:stepLeray} degenerates at the $E_2$-page.

By construction,
$\pi_{\bi_k}^* \mathcal F(M) \cong \cE_{\bi_k}(M)$, and hence
\[
H^0\bigl( Z(\bi_k), \pi_{\bi_k}^* \mathcal F(M) \bigr) \cong H^0\bigl( Z(\bi_k), \cE_{\bi_k}(M) \bigr) \cong \mathscr{C}_{\hs_k(\Psi), k}(M)^{\vee}.
\]

We thus conclude that
\[
H^0\bigl( \tbP(k) / (\tbP(k) \cap \tbP(k+1)),  \mathcal F(M) \bigr) \cong \mathscr{C}_{\hs_k(\Psi), k}(M)^{\vee},
\]
which is the dual form of the asserted identity.
\end{proof}

\begin{prop}\label{prop:idY(k)}
Let $\Psi \subset \Delta^+$ be a root ideal, and let $e_1(\Psi) \le k < n$.
Then the variety $\sY_\Psi(k)$ admits the structure of a $\tbP(k)$-equivariant $\sY_\Psi(k+1)$-fibration over the base $\tbP(k) / ( \tbP(k) \cap \tbP(k+1) )$.
\end{prop}

\begin{proof}
Set $\bQ ( k ) :=\tbP(k)\cap\tbP(k+1)$. By repeated applications of Corollary~\ref{cor:D-branch}, there is a natural surjection
\begin{equation}
\bigoplus_{\la \in \Par} N^{\Psi}(\la; k)^{\vee} \longrightarrow \bigoplus_{\la \in \Par} N^{\Psi}(\la; k+1)^{\vee} \otimes \C_{-m_k(\la)\La_k}.\label{eqn:N(k)-surj}
\end{equation}
Lemma~\ref{lem:ind}, applied degree by degree, identifies
the left-hand side (including its multiplication) with the global sections
of the $\Par$-graded sheaf $\mathcal A_k$ of algebras on $\tbP(k) / \bQ ( k )$ given as the sheaf of sections of the direct sum of vector bundles
\[
\bigoplus_{\la\in\Par} \tbP(k)\times^{\bQ ( k )}
 \left(N^\Psi(\la;k+1)
       \otimes\C_{m_k(\la)\La_k}\right) \rightarrow \tbP(k) / \bQ ( k ).
\]
Consequently,
\[
 \Proj_{\tbP(k) / \bQ ( k ),\Par}\mathcal A_k
 \cong\tbP(k)\times^{\bQ ( k )}\sY_\Psi(k+1).
\]
The evaluation homomorphism of the section algebra therefore induces a
$\tbP(k)$-equivariant morphism
\begin{equation}
\pi_k \colon \tbP(k) \times^{\bQ ( k )} \sY_\Psi(k+1) \longrightarrow \sY_\Psi(k).\label{eqn:defpik}
\end{equation}
The surjection~\eqref{eqn:N(k)-surj} of graded algebras also identifies the fiber over
the identity coset $\bQ(k)/\bQ(k)$ with a closed subscheme $\sY_\Psi(k+1)\subset
\sY_\Psi(k)$ through $\pi_k$.

To analyze the image of $\sY_\Psi(k+1)$, consider weights of the form $c \varpi_k$ for $c \in \Z_{\ge 0}$.  
Since $m_j(c \varpi_k) = 0$ for all $j > k$,  
repeated applications of Corollary~\ref{cor:DLcomm} yield
\[
N^\Psi(c \varpi_k; k+1) = \C \qquad \text{for all } c \in \Z_{\ge 0}.
\]

Consider the $\tbP(k)$-equivariant morphism
\[
\psi \colon \sY_\Psi(k) \longrightarrow \P(L(\La_k)),
\]
induced by the line bundle $\mathcal{O}_{\sY_\Psi(k)}(\varpi_k)$.  On the
closed subscheme $\sY_\Psi(k+1)$, its degree-$c$ section space is the
one-dimensional space $\C_{-c\La_k}$.  Hence the whole subscheme is
mapped scheme-theoretically to the point
\begin{equation}
\mathrm{pt} = \Proj_{\Z_{\ge 0}} \bigoplus_{c \ge 0} \C_{-c \La_k} \subset \P(N^\Psi(\varpi_k; k)),\label{eqn:ptp}
\end{equation}
which is fixed by $\bQ ( k )$. Under the inclusion
$N^\Psi(\varpi_k;k)\subset L(\La_k)$ of
Theorem~\ref{thm:DCF}(2), the generator of
$\C_{-\La_k}^\vee$ is the highest-weight vector
$\bv_{\La_k}$. Indeed, under the identification
$G(k)\cong\SL(\hs_k(\Psi)+n-k+1)$, the line
$\C\bv_{\La_k}$ is the highest-weight line in the corresponding
fundamental representation. Its stabilizer is the maximal parabolic
$G(k)\cap\tbP(k+1)$ appearing in Lemma~\ref{lem:Ustr}.
Since $\tbP(k)=G(k)\tbB$, its stabilizer in $\tbP(k)$ is $\bQ(k)$. Thus the point is $[\bv_{\La_k}]$, and its stabilizer in
$\tbP(k)$ is exactly $\bQ ( k )$.

It follows that the $\tbP(k)$-orbit of this point is
\[
\tbP(k)/\bQ ( k ) \cong \P^{\hs_k(\Psi) + n - k} = \P(N^\Psi(\varpi_k; k)).
\]
The morphism $\psi$ therefore factors through a $\tbP(k)$-equivariant
morphism $\sY_\Psi(k)\to \tbP(k) / \bQ ( k )$, and its fiber over the identity
coset is $\sY_\Psi(k+1)$ by the preceding graded-algebra calculation.

Thus the $\bQ(k)$-equivariant identification of this fiber with
$\sY_\Psi(k+1)$ induces an isomorphism
\[
\tbP(k)\times^{\bQ(k)}\sY_\Psi(k+1)
\stackrel{\cong}{\longrightarrow}
\sY_\Psi(k),
\]
which coincides with $\pi_k$, as required.
\end{proof}

\begin{cor}\label{cor:Y(k)closed}
Let $\Psi \subset \Delta^+$ be a root ideal, and let $e_1(\Psi) \le k \le n$. 
Then the map
\begin{equation}
\sY_\Psi(k) \longrightarrow \prod_{j=k}^n \P \bigl( L ( \La_j ) \bigr),\label{eqn:partialX}
\end{equation}
induced from Corollary~\ref{cor:XmapP} and Lemma~\ref{lem:idXsX}, is a closed immersion.
\end{cor}

\begin{proof}
We argue by downward induction on $k$.  For $k=n$, the definition and
Lemma~\ref{lem:ind} identify $\sY_\Psi(n)$ with
\[
\tbP(n)/(\tbP(n)\cap\tbP(n+1))
\]
and~\eqref{eqn:partialX} with its closed-orbit embedding into
$\P(L(\La_n))$.

Assume the assertion holds for $k+1$.  By
Proposition~\ref{prop:idY(k)} and its proof, $\sY_\Psi(k)$ is a
$\tbP(k)$-equivariant $\sY_\Psi(k+1)$-fibration over
\[
\tbP(k)/(\tbP(k)\cap\tbP(k+1)),
\]
and the composite of~\eqref{eqn:partialX} with the projection onto
$\P(L(\La_k))$ is the fibration projection followed by the
closed-orbit embedding. Its fiber over $[\bv_{\La_k}]$ is
$\sY_\Psi(k+1)$, on which~\eqref{eqn:partialX} is a closed immersion
by the induction hypothesis.

The associated-bundle description and the $\tbP(k)$-equivariance
therefore identify~\eqref{eqn:partialX}, over each local
trivialization and up to an automorphism induced by the $\tbP(k)$-action
on $\prod_{j=k+1}^n\P(L(\La_j))$, with the product of the identity on
the base and this closed immersion of the fiber.  Hence
\eqref{eqn:partialX} is a closed immersion.
\end{proof}

\begin{cor}\label{cor:Y(k)sm}
Let $\Psi \subset \Delta^+$ be a root ideal, and let $e_1(\Psi) \le k \le n$. Then the variety $\sY_\Psi(k)$ is smooth, and the dimension of $\sX_\Psi(w_0^\Psi)$ is given by
\begin{equation}
\dim \sX_\Psi(w_0^\Psi) = |\Psi| + \ell(w_0^\Psi).
\end{equation}
\end{cor}

\begin{proof}
By Lemma~\ref{lem:ind} and Proposition~\ref{prop:idY(k)}, the variety $\sY_\Psi(k)$ is a successive projective-space bundle of dimensions $\{\hs_j(\Psi) + n - j\}_{j = k}^n$, and is therefore smooth.

To compute the dimension, we observe that
\begin{align*}
\dim \sX_\Psi(w_0^\Psi) &= \dim \sY_\Psi(e_1(\Psi)) = \sum_{k = e_1(\Psi)}^n (\hs_k(\Psi) + n - k) \\
&= \sum_{k = e_1(\Psi)}^n \hs_k(\Psi) + \sum_{i = 0}^{n - d_1(\Psi) - 1} i \\
&= |\Psi| + \frac{(n - d_1(\Psi))(n - d_1(\Psi) - 1)}{2} = |\Psi| + \ell(w_0^\Psi).
\end{align*}
Here we apply~\eqref{eqn:Psicard} to evaluate the sum $\sum \hs_k(\Psi)$, and recall that $w_0^\Psi$ is the longest element in $\Sym_{n - d_1(\Psi)}$. This completes the proof.
\end{proof}

\begin{thm}\label{thm:str}
Let $\Psi \subset \Delta^+$ be a root ideal. Then there exists a $G$-equivariant closed immersion
\begin{equation}
\sX_{\Psi}(w_0) \cong G \times^Q \sX_\Psi(w_0^\Psi) \hookrightarrow \prod_{i=1}^n \P\big( L(\La_i) \big),\label{eqn:Xwemb}
\end{equation}
where $Q \subset G$ denotes the parabolic subgroup generated by $B$ and $P_i$ for each $e_1(\Psi) \le i < n$.

Moreover, for $\Psi$-tame $w\in\Sym_n$, we have
\begin{equation}\label{eqn:Xw-associated}
\sX_\Psi(w)
\cong
X_Q(w)\times_{G/Q} \sX_\Psi(w_0) \quad \text{for} \quad X_Q(w) := \overline{BwQ/Q} \subset G/Q.
\end{equation}
As a consequence,
\[
\dim\sX_\Psi(w)=\ell(w)+|\Psi|,
\]
and $\sX_\Psi(w)$ is smooth if and only if the Schubert variety $X(w)$
is smooth.
\end{thm}

\begin{proof}
For $w\in\Sym_n$, put
\[
R_\Psi(w):=
\bigoplus_{\la\in\Par}N^\Psi_w(\la)^\vee.
\]
By Lemma~\ref{lem:tame}, the $\Par$-graded algebra
$R_\Psi(w_0^\Psi)$ is stable under $Q$. Let $\mathcal R_Q$ be the
$\Par$-graded sheaf of algebras on $G/Q$ given as the sheaf of sections
of the direct sum of vector bundles
\[
\bigoplus_{\la\in\Par}
G\times^Q N^\Psi_{w_0^\Psi}(\la)
\longrightarrow G/Q.
\]
Its relative multigraded Proj is canonically
\begin{equation}\label{eqn:relative-proj-Q}
\Proj_{G/Q,\Par}\mathcal R_Q
\cong G\times^Q \sX_\Psi(w_0^\Psi).
\end{equation}

Let $\bi_0$ record a reduced expression of $w_0$, and denote by
\[
\eta_0\colon Z(\bi_0)\longrightarrow G/B\longrightarrow G/Q
\]
the natural morphism. We have
\[
(\eta_0)_*\cO_{Z(\bi_0)}=\cO_{G/Q},
\qquad \text{and} \qquad
\mathbb R^{>0}(\eta_0)_*\cO_{Z(\bi_0)}=0
\]
(cf.~\cite[Theorem~8.2.2(c) and A.24]{Kum02}).
The geometric realization of Demazure functors, together with the projection formula, gives an isomorphism of
$\Par$-graded algebras
\begin{align}
R_\Psi(w_0)
&\cong \mathscr D_{\bi_0}^\dag ( R_\Psi(w_0^\Psi) ) = 
H^0\bigl(Z(\bi_0),\eta_0^*\mathcal R_Q\bigr)
\cong
H^0\bigl(G/Q,\mathcal R_Q\bigr).
\label{eqn:Rw0-associated}
\end{align}

We next describe the morphism defined by the fundamental components of
$\mathcal R_Q$. For $e_1(\Psi)\le i\le n$, compose the morphism $f_i$
in~\eqref{eqn:deffk} with the linear closed immersion
\[
\P\bigl(N^\Psi_{w_0^\Psi}(\varpi_i)\bigr)
\hookrightarrow\P(L(\La_i)),
\]
and continue to denote the resulting morphism by $f_i$. By the
$Q$-action on $R_\Psi(w_0^\Psi)$, each $f_i$ is $Q$-equivariant.
For $1\le i \le d_1(\Psi)$, we have
\[
N^\Psi_{w_0^\Psi}(\varpi_i)=\C\bv_{\La_i}
\]
by Theorem~\ref{thm:DCF}(2). We therefore obtain a $G$-equivariant
morphism
\begin{equation}\label{eqn:associated-embedding}
\iota\colon
G\times^Q \sX_\Psi(w_0^\Psi)
\longrightarrow
\prod_{i=1}^n\P(L(\La_i))
\end{equation}
whose components are
\[
[g,x]\longmapsto
\begin{cases}
[g\bv_{\La_i}] & (1\le i\le d_1(\Psi)),\\
g f_i(x) & (e_1(\Psi)\le i\le n).
\end{cases}
\]
The $G$-stabilizer of $\bigl([\bv_{\La_i}]\bigr)_{i=1}^{d_1(\Psi)}$ is $Q$. Hence the first $d_1(\Psi)$ components identify $G/Q$ with a
closed $G$-orbit.

If $e_1(\Psi) = n+1$, then we have $\sX_\Psi(w_0^\Psi)=\mathrm{pt}$ and~\eqref{eqn:associated-embedding} yields the closed immersion~\eqref{eqn:Xwemb}.

If $e_1(\Psi) \le n$, then the remaining components of~\eqref{eqn:associated-embedding} restrict on the fiber over the identity coset to a closed immersion
\[
\sX_\Psi(w_0^\Psi)=\sY_\Psi(e_1(\Psi))
\longrightarrow
\prod_{i=e_1(\Psi)}^n\P(L(\La_i))
\]
by Corollary~\ref{cor:Y(k)closed}. Over every local trivialization of $G\times^Q \sX_\Psi(w_0^\Psi) \to G/Q$, the morphism
$\iota$ is therefore the product of the identity on the base with this
closed immersion, up to an automorphism of
$\prod_{i=1}^n\P(L(\La_i))$
induced by the diagonal $G$-action. Hence $\iota$ is a closed immersion.

The line bundles defining $\iota$ are the fundamental components of
$\mathcal R_Q$. Thus~\eqref{eqn:Rw0-associated} identifies their
multisection ring with $R_\Psi(w_0)$. Since they define the closed
immersion~\eqref{eqn:associated-embedding}, the multigraded Proj
construction yields
\[
\sX_\Psi(w_0)
=\Proj_\Par R_\Psi(w_0)
\cong
G\times^Q \sX_\Psi(w_0^\Psi).
\]
Under this isomorphism, $\iota$ gives the closed immersion
\eqref{eqn:Xwemb}.

We denote the resulting bundle projection by
\[
\mathsf{pr}\colon
\sX_\Psi(w_0)\cong G\times^Q \sX_\Psi(w_0^\Psi)
\longrightarrow G/Q,
\]
which is an identity when $e_1(\Psi) = n+1$.

We now consider a $\Psi$-tame element $w$. Write
\[
w=vw_0^\Psi,
\qquad
\ell(w)=\ell(v)+\ell(w_0^\Psi).
\]
We have $X_Q (w) = \overline{BvQ/Q}$ and $\dim X_Q(w) = \ell (v)$. Choose a reduced expression of $v$, let $\bi'$ be the corresponding
sequence, and let
\[
h\colon Z(\bi')\longrightarrow X_Q(w)
\]
be the Bott--Samelson resolution. Set
\[
\mathcal R_{Q,w}:=\mathcal R_Q|_{X_Q(w)}.
\]
Since $w=vw_0^\Psi$ and $R_\Psi(w_0^\Psi)$ is $Q$-stable,
the definition of $N^\Psi_w(\la)$ gives
\[
R_\Psi(w)
\cong
\mathscr D_v^\dag\bigl(R_\Psi(w_0^\Psi)\bigr).
\]
The geometric realization of the Demazure functors and the projection
formula for $h$ give, compatibly with multiplication,
\begin{equation}\label{eqn:Rw-relative}
R_\Psi(w)
\cong
H^0\bigl(Z(\bi'),h^*\mathcal R_{Q,w}\bigr)
\cong
H^0\bigl(X_Q(w),\mathcal R_{Q,w}\bigr),
\end{equation}
where we use
\[
h_*\cO_{Z(\bi')}=\cO_{X_Q(w)} \qquad \text{and} \qquad \mathbb R^{>0} h_*\cO_{Z(\bi')}=0
\]
(cf.~\cite[Theorem~8.2.2(c) and A.24]{Kum02}).

On the other hand, relative multigraded Proj commutes with restriction
to $X_Q(w)$, and hence
\[
\Proj_{X_Q(w),\Par}\mathcal R_{Q,w}
\cong
X_Q(w)\times_{G/Q}
\Proj_{G/Q,\Par}\mathcal R_Q\cong
\mathsf{pr}^{-1}\bigl(X_Q(w)\bigr).
\]
The restriction of~\eqref{eqn:associated-embedding} to this relative
Proj is again a closed immersion, and~\eqref{eqn:Rw-relative}
identifies its multisection ring with $R_\Psi(w)$. Therefore
\[
\sX_\Psi(w)
\cong
\mathsf{pr}^{-1}\bigl(X_Q(w)\bigr)
\]
scheme-theoretically, that is~\eqref{eqn:Xw-associated}.

The restriction of $\mathsf{pr}$ to
$\mathsf{pr}^{-1}(X_Q(w))$ is a locally trivial fibration with fiber $\sX_\Psi(w_0^\Psi)$. By Corollary~\ref{cor:Y(k)sm} (when $e_1(\Psi) \le n$) or by direct inspection (when $e_1(\Psi) = n+1$),
\[
\dim \sX_\Psi(w_0^\Psi)=|\Psi|+\ell(w_0^\Psi),
\]
and consequently
\[
\dim\sX_\Psi(w)
=
\dim X_Q(w)+\dim \sX_\Psi(w_0^\Psi)
=
\ell(v)+\ell(w_0^\Psi)+|\Psi|
=
\ell(w)+|\Psi|.
\]

The fiber $\sX_\Psi(w_0^\Psi)$ is smooth, so the locally trivial fibration
\eqref{eqn:Xw-associated} shows that $\sX_\Psi(w)$ is smooth if and
only if $X_Q(w)$ is smooth. Since $w=vw_0^\Psi$ is the longest element
in its coset modulo the Weyl group of $Q$, the projection
$G/B\to G/Q$ restricts to a locally trivial fibration
\[
X(w)\longrightarrow X_Q(w)
\]
with fiber $Q/B$. Thus $X_Q(w)$ is smooth if and only if $X(w)$ is
smooth, completing the proof.
\end{proof}

\begin{ex}[$n=4$]\label{ex:n4}
We illustrate the construction of $\sX_\Psi(w_0)$ in the case $G = \GL(4,\C)$, using the root ideal
\[
\Psi = \{\epsilon_1 - \epsilon_3, \epsilon_1 - \epsilon_4,\epsilon_2 - \epsilon_4\}.
\]
In this case, we have $e_1(\Psi) = 3$ and
\[
\hs_2(\Psi) = 0, \qquad \hs_3(\Psi) = 1, \qquad \hs_4(\Psi) = 2.
\]

We begin by defining the subspace
\[
V^{(4)} := \C \bv_{1111} \oplus \C \bv_{2110} \oplus \C \bv_{1210} \subset L(\La_4),
\]
where $\bv_{1111} := \bv_{\La_4}$ is the highest-weight vector, and the remaining vectors $\bv_{2110}$ and $\bv_{1210}$ have $\tT$-weights $\La_4 - \alpha_0$ and $\La_4 - \alpha_0 - \alpha_1$, respectively. The vector $\bv_{1111}$ has degree~$0$, while $\bv_{2110}$ and $\bv_{1210}$ have degree~$-1$. We then have
\[
\P^2 \cong \sY_\Psi(4) = \P(V^{(4)}) \subset \P(L(\La_4)),
\]
which is preserved under
\[
\tbP(4) := \left< \SL(2,0),\, \SL(2,1),\, \tbB \right> \subset \widetilde{G}(\!(z)\!).
\]

Next, set
\[
V^{(3)} := \C \bv_{1110} \oplus \C \bv_{1101} \oplus \C \bv_{2100} \subset L(\La_3),
\]
where $\bv_{1110} := \bv_{\La_3}$ is the highest-weight vector, and the remaining vectors $\bv_{1101}$ and $\bv_{2100}$ have $\tT$-weights $\La_3 - \alpha_3$ and $\La_3 - \alpha_3 - \alpha_0$, respectively. Here, $\bv_{1110}$ and $\bv_{1101}$ have degree~$0$, while $\bv_{2100}$ has degree~$-1$. The projective space $\P(V^{(3)}) \subset \P(L(\La_3))$ is preserved under
\[
\tbP(3) := \left< \SL(2,3), \SL(2,0), \tbB \right> \subset \widetilde{G}(\!(z)\!).
\]
Define $
G(3) := \left< \SL(2,3), \SL(2,0) \right> \cong \SL(3) \subset \widetilde{G}(\!(z)\!)$ and let $P(3) \subset G(3)$ be the parabolic subgroup stabilizing $V^{(4)}$. Then $\sY_\Psi(3)$ admits the structure
\[
\sY_\Psi(3) = G(3) \times^{P(3)} \P(V^{(4)}) \cong G(3) \cdot ([\bv_{1110}] \times \P(V^{(4)})) \subset \P(L(\La_3)) \times \P(L(\La_4)),
\]
which defines a $\P(V^{(4)})$-bundle over $\P(V^{(3)})$.

To describe the $G(3)$-orbit $G(3) \cdot \P(V^{(4)})$, we extend $V^{(4)}$ to include three additional vectors
\[
\bv_{2101}, \bv_{1201} \text{ (degree $-1$), and } \bv_{2200} \text{ (degree $-2$)},
\]
determined by their $\tT$-weights. This results in
\[
W^{(4)} := V^{(4)} \oplus \C \bv_{2101} \oplus \C \bv_{1201} \oplus \C \bv_{2200},
\]
which is stable under the action of $G(3)$ and $\tbB$.

We then have the following projective embedding:
\begin{align*}
\sY_\Psi(3) & = \Bigg\{ \left( \begin{matrix} x_{1110}^{(3)} \\ x_{1101}^{(3)} \\ {\color{red}x_{2100}^{(3)}} \end{matrix} \right) 
\parallel
\left( \begin{matrix} {\color{red}x_{1210}^{(4)}} \\ {\color{red}x_{1201}^{(4)}} \\ {\color{blue}x_{2200}^{(4)}} \end{matrix} \right),
x_{1110}^{(3)} {\color{red}x_{2101}^{(4)}} + x_{1101}^{(3)} {\color{red}x_{2110}^{(4)}} + {\color{red}x_{2100}^{(3)}} x_{1111}^{(4)} = 0 \Bigg\} \\
& \hskip 20mm \subset \Big\{ ([x_{\bullet}^{(3)}], [x_{\bullet}^{(4)}]) \in \P(V^{(3)}) \times \P(W^{(4)}) \Big\}  \cong \P^2 \times \P^5,
\end{align*}
where $x_\bullet^{(i)}$ denotes the homogeneous coordinate dual to
$\bv_\bullet^{(i)}$, and the coloring indicates the degree of the corresponding coordinate function: black for~$0$, red for~$1$, and blue for~$2$.

The degree-zero locus is the subvariety
\[
\P^1 \cong \P(\C \bv_{1110} \oplus \C \bv_{1101}) \times \P(\C \bv_{1111}) \subset \P(V^{(3)}) \times \P(W^{(4)}),
\]
noting that $\P(\C \bv_{1111})$ is a point.

We now describe the attracting locally closed subset $\mathcal U^-$ of the $\tbB$-fixed point $([\bv_{1110}], [\bv_{1111}])$ by setting
\[
x_{1101}^{(3)} = 0, \qquad x_{1110}^{(3)} = 1 = x_{1111}^{(4)}.
\]
The coordinates ${\color{red}x_{2110}^{(4)}}, {\color{red}x_{1210}^{(4)}}, {\color{red}x_{2101}^{(4)}}$ are then free, while the others are determined by
\[
{\color{red}x_{1201}^{(4)}} = \frac{x_{1101}^{(3)} {\color{red}x_{1210}^{(4)}}}{x_{1110}^{(3)}}= 0, \quad {\color{blue}x_{2200}^{(4)}} = \frac{{\color{red}x_{1210}^{(4)}}{\color{red}x_{2100}^{(3)}}}{x_{1110}^{(3)}}, \quad  {\color{red}x_{2100}^{(3)}} = - \frac{x_{1110}^{(3)} \color{red}x_{2101}^{(4)}}{x_{1111}^{(4)}} - \frac{x_{1101}^{(3)} {\color{red}x_{2110}^{(4)}}}{x_{1111}^{(4)}}.
\]
These variables have $T$-weights $-\epsilon_1 + \epsilon_4, -\epsilon_2 + \epsilon_4, -\epsilon_1 + \epsilon_3$, respectively. We thus obtain a $B$-equivariant (degree-preserving) identification
\[
\mathcal U^- = \exp \left( \C E_{14} z^{-1} + \C E_{24} z^{-1} + \C E_{13} z^{-1} \right) \cdot ([\bv_{1110}], [\bv_{1111}]) \subset \P(V^{(3)}) \times \P(W^{(4)}).
\]

Since $\sY_\Psi(3)$ is $\SL(2,3)$-stable, we conclude that
\begin{equation}
T^*_\Psi X \cong G \times^B \mathcal U^- \subset G \times^{P_3} \sY_\Psi(3) \cong \sX_\Psi(w_0).\label{eqn:exinj}
\end{equation}

Finally, the relation
\[
x_{1110}^{(3)} {\color{red}x_{2101}^{(4)}} + x_{1101}^{(3)} {\color{red}x_{2110}^{(4)}} + {\color{red}x_{2100}^{(3)}} x_{1111}^{(4)} = 0,
\]
shows that the complement $\sX_\Psi(w_0) \setminus T^*_\Psi X$ is given by the locus $x_{1111}^{(4)} = 0$.
\end{ex}

\section{Properties of the variety $\sX_\Psi$}

We continue to work within the framework established in the previous section.

\begin{thm}\label{thm:Xmain}
Let $\Psi \subset \Delta^+$ be a root ideal, and let $w \in \Sym_n$ be $\Psi$-tame. For each $\la \in \Par$, we have:
\begin{enumerate}
\item $H^{>0} \bigl( \sX_{\Psi} ( w ), \cO_{\sX_\Psi ( w )} ( \la )\bigr) = 0$;
\item there is an isomorphism of $\tbB$-modules
\[
H^0\bigl(\sX_\Psi(w),\cO_{\sX_\Psi(w)}(\la)\bigr)^\vee
\cong N_w^\Psi(\la);
\]
\item if $\la_1 > 0$, then the module $N_w^\Psi(\la)$ admits a
$D^{(\la_1)}$-filtration. In all cases,
$H^0(\sX_\Psi(w),\cO_{\sX_\Psi(w)}(\la))$, regarded as a
$B$-module, admits an excellent filtration in the sense of
van~der~Kallen~\cite{vdK89}.
\end{enumerate}
\end{thm}

\begin{proof}
Choose $v \in \Sym_n$ such that $w=vw_0^\Psi$ and $\ell(w)=\ell(v)+\ell(w_0^\Psi)$.
Repeated applications of Lemma~\ref{lem:tame} give
\[
N_w^\Psi(\la)\cong N_v^\Psi(\la)
\qquad (\la\in\Par).
\]
Let $\bi=\bi(\Psi,v)$ be the sequence introduced at~\eqref{eqn:defbiPw} in Section~\ref{sec:rot}. Thus $\bi$ records, in order, a reduced expression
of $v$ and the minimal parabolics occurring in the definition of
$N_v^\Psi(\la)$. Let $Z(\bi)$ be the corresponding Bott--Samelson
variety.

By~\eqref{eqn:section-ring-BS}, we have an isomorphism of
$\Par$-graded algebras
\begin{equation}\label{eqn:section-ring-Xmain}
\bigoplus_{\la\in\Par}N_w^\Psi(\la)^\vee
\cong
\bigoplus_{\la\in\Par}
H^0\bigl(Z(\bi),\mathcal L_v^\Psi(\la)\bigr).
\end{equation}
The fundamental components of this multisection ring are base-point-free by Theorem~\ref{thm:DCF}(3).
Their evaluation maps therefore define a morphism
\begin{equation}\label{eqn:pi-Xmain}
\pi_{\bi}\colon Z(\bi)\longrightarrow\sX_\Psi(w)
\end{equation}
such that
\begin{equation}\label{eqn:pi-pullback}
\pi_{\bi}^*\cO_{\sX_\Psi(w)}(\la)
\cong\mathcal L_v^\Psi(\la)
\qquad (\la\in\Par).
\end{equation}
Here the right-hand side of~\eqref{eqn:pi-Xmain} is precisely the multigraded
$\Proj$ of the algebra on the left-hand side of
\eqref{eqn:section-ring-Xmain}.  Under~\eqref{eqn:Xwemb} and~\eqref{eqn:Xprojemb}, the morphism $\pi_{\bi}$ is precisely the one
defined by the fundamental evaluation maps.

The associated-bundle descriptions in
Proposition~\ref{prop:idY(k)} and Theorem~\ref{thm:str} show that the sequence $\bi$ is obtained by
concatenating the sequence defining the usual Bott--Samelson resolution
of $X_Q(w)$ with the sequences defining the Bott--Samelson resolutions
of the projective homogeneous spaces occurring as the successive
fibers. Consequently, $\pi_{\bi}$ is projective and birational.
Since $Z(\bi)$ is smooth, $\pi_{\bi}$ is a resolution of singularities.

Moreover, Schubert varieties have
rational singularities~\cite[Theorem~8.2.2(c)]{Kum02}, and
\eqref{eqn:Xw-associated} identifies
$\sX_\Psi(w)$ locally over $X_Q(w)$ with the product of $X_Q(w)$ and the
smooth variety $\sX_\Psi(w_0^\Psi)$. Hence $\sX_\Psi(w)$ has rational
singularities. It follows from~\cite[Definition~5.8 and Theorem~5.10]{KM98} that
\begin{equation}\label{eqn:RpiO}
(\pi_{\bi})_*\cO_{Z(\bi)}
\cong\cO_{\sX_\Psi(w)} \quad \text{and} \quad \mathbb R^{>0}(\pi_{\bi})_*\cO_{Z(\bi)} = 0.
\end{equation}

The projection formula,~\eqref{eqn:pi-pullback}, and
Proposition~\ref{prop:inclD} now give
\[
H^q\bigl(\sX_\Psi(w),\cO_{\sX_\Psi(w)}(\la)\bigr)
\cong
H^q\bigl(Z(\bi),\mathcal L_v^\Psi(\la)\bigr)
\cong
\begin{cases}
N_w^\Psi(\la)^\vee & (q=0),\\
0 & (q>0).
\end{cases}
\]
This proves (1) and (2).

If $\la_1=0$, then $\la=0$ and
\[
 H^0(\sX_\Psi(w),\cO_{\sX_\Psi(w)})\cong\C,
\]
so the excellent-filtration assertion is immediate. Assume that
$\la_1>0$. Repeated applications of Corollary~\ref{cor:D-branch} to the
definition of $N_w^\Psi(\la)$ give a $D^{(\la_1)}$-filtration. By Theorem~\ref{thm:D-branch}, for every $\mu\in\sP$, $D_\mu^{(k)}\otimes\C_\wp$ admits a $D^{(k+1)}$-filtration.
Since $\C_\wp$ is trivial on $B$, the same holds for
$D_\mu^{(k)}$ as a $B$-module. For sufficiently large
$k'$, every $D^{(k')}_\mu$ occurring in this process is a Demazure module
for $G$; see~\cite[\S3.5]{Jos85},~\cite[Theorem~1]{FL07}, or~\cite[Theorem~B]{Kat26b}. Dualizing the resulting filtration gives an
excellent filtration of
$H^0\bigl(\sX_\Psi(w),\cO_{\sX_\Psi(w)}(\la)\bigr)$, proving (3).
\end{proof}

\begin{cor}\label{cor:BMPid}
Let $\Psi \subset \Delta^+$ be a root ideal, and let $w \in \Sym_n$ be $\Psi$-tame. For each $\la \in \Par$, we have
\[
\gch H^0 \bigl( \sX_{\Psi} ( w ), \cO_{\sX_\Psi ( w )} ( \la ) \bigr)^{\vee} = \bigl[ H ( \Psi; \la; w ) \bigr]_{q \mapsto q^{-1}}.
\]
\end{cor}

\begin{proof}
This follows immediately by combining Theorem~\ref{thm:Xmain} with Theorem~\ref{thm:BMP}.
\end{proof}

For $1\le i\le n$, let
\[
V(\varpi_i) \cong L ( \La_i )_0  \subset L(\La_i)
\]
be the degree-zero subspace of $L(\La_i)$.
Note that $V(\varpi_n)  = \C\bv_{\La_n}$ is the determinant character of $G$.
Thus we have a $G$-equivariant, $\Gm$-stable closed immersion
\begin{equation}\label{eqn:finite-fixed-product}
\prod_{i=1}^n\P(V(\varpi_i) )
\hookrightarrow
\prod_{i=1}^n\P(L(\La_i)),
\end{equation}
where the last factor on the left is a point.

Since every
$L(\La_i)$ is supported in degrees at most zero, the attracting set
of $\P(V (\varpi_i) )$ in $\P(L(\La_i))$ for $t\to\infty$ is open.

\begin{lem}\label{lem:Xtop}
Let $\Psi \subset \Delta^+$ be a nonempty root ideal, and let $w \in \Sym_n$ be
$\Psi$-tame. Under the embedding induced by~\eqref{eqn:Xwemb}, we have
\[
\sX_\Psi(w)\cap\prod_{i=1}^n\P( L ( \La_i )_0 )
\cong X(w).
\]
\end{lem}

\begin{proof}
For $1\le j\le d_1(\Psi)$, Theorem~\ref{thm:DCF}(2) gives
\[
N_{w_0^\Psi}^\Psi(\varpi_j)=\C\bv_{\La_j}.
\]
Consequently, the $j$-th projection of $\sX_\Psi(w_0^{\Psi})$ is the single point
$[\bv_{\La_j}]$. Under the embedding
$V(\varpi_j)\hookrightarrow L(\La_j)$, this is precisely the image
of the highest-weight point $[\bv_{\varpi_j}]$. Therefore the intersection of the $j$-th image of $\sX_\Psi(w_0^{\Psi})$ with
$\P(V(\varpi_j))$ consists of that point.

For $e_1(\Psi)\le k\le n$, let $f_{\le k}$ denote the product of the
components $f_j$ in~\eqref{eqn:deffk}, followed by their linear
inclusions into $\P(L(\La_j))$, for $e_1(\Psi)\le j\le k$. Put
\[
K(k):=\left\langle\SL(2,i)\mid k\le i<n\right\rangle
\cong\SL(n-k+1),
\]
with $K(n)=\{1\}$. We claim that
\begin{equation}\label{eqn:fixed-partial-image}
\operatorname{Im}f_{\le k}
\cap
\prod_{j=e_1(\Psi)}^k\P(L(\La_j)_0)
=
K(e_1(\Psi))\cdot
\bigl([\bv_{\varpi_j}]\bigr)_{j=e_1(\Psi)}^k.
\end{equation}

For $k=e_1(\Psi)<n$, Proposition~\ref{prop:idY(k)} identifies
$\operatorname{Im}f_{\le e_1(\Psi)}$ with the projectivization of the
fundamental representation of $G(e_1(\Psi))$. When
$e_1(\Psi)=n$, the same identification follows from
Lemma~\ref{lem:ind}. Its degree-zero subspace is the corresponding
fundamental representation of $K(e_1(\Psi))$. Hence its intersection
with $\P(V( \varpi_{e_1(\Psi)} ) )$ is the
$K(e_1(\Psi))$-orbit of the highest-weight point.

Suppose~\eqref{eqn:fixed-partial-image} holds for some $k<n$. Over the
highest-weight point in its right-hand side, the next stage is the
projective homogeneous space supplied by Proposition~\ref{prop:idY(k)}.
Its intersection with $\P(V( \varpi_{k+1} ))$ is the orbit of its
highest-weight point under
\[
K(k+1)\cong\SL(n-k).
\]
Equivariance then proves~\eqref{eqn:fixed-partial-image} for $k+1$.
This completes the induction.

Taking $k=n$ and using Corollary~\ref{cor:Y(k)closed}, we obtain
\begin{equation}\label{eqn:F-fixed}
\sX_\Psi(w_0^\Psi) \cap\prod_{i=e_1(\Psi)}^n\P(L(\La_i)_0)
=X(w_0^\Psi).
\end{equation}
Indeed, the orbit on the right-hand side of
\eqref{eqn:fixed-partial-image} is the full flag variety of
$K(e_1(\Psi))\cong\SL(n-d_1(\Psi))$.

By Theorem~\ref{thm:str},
\[
\sX_\Psi(w_0)\cong G\times^Q\sX_\Psi(w_0^\Psi).
\]
Combining this with~\eqref{eqn:F-fixed}, we obtain
\begin{equation}
\sX_\Psi(w_0)\cap\prod_{i=1}^n\P(L (\La_i)_0)
\cong
G\times^QX(w_0^\Psi)
\cong G/B.\label{eqn:deg0int}
\end{equation}
Finally,~\eqref{eqn:Xw-associated} gives
\[
\sX_\Psi(w)
\cong
X_Q(w)\times_{G/Q}\sX_\Psi(w_0).
\]
Intersecting with~\eqref{eqn:deg0int}, we therefore obtain
\[
\sX_\Psi(w)\cap\prod_{i=1}^n\P(L (\La_i)_0)
\cong
X_Q(w)\times_{G/Q}G/B
=X(w). \qedhere
\]
\end{proof}

\begin{thm}\label{thm:XTid}
For every root ideal $\Psi\subset\Delta^+$, the $\Gm$-attracting set of
the fixed component
\[
X=X(w_0)\subset\sX_\Psi(w_0)
\]
is open and dense, and is $G$-equivariantly isomorphic to
$T^*_\Psi X$.
\end{thm}

\begin{proof}
If $\Psi=\varnothing$, then $e_1(\Psi)=n+1$ and
$\sX_\Psi(w_0^\Psi)=\mathrm{pt}$. Hence Theorem~\ref{thm:str} gives
$\sX_\Psi(w_0)=X$, and the assertion is immediate.

Assume from now
on that $\Psi\ne\varnothing$.

For $1\le i\le n$, let $\mathfrak A_i$ be the attracting set of
$\P( L(\La_i)_0 )$ in $\P(L(\La_i))$ for $t\to\infty$. Since all
nonzero degrees of $L(\La_i)$ are negative, $\mathfrak A_i$ is open.
Lemma~\ref{lem:Xtop} and the $\Gm$-equivariant closed immersion
\eqref{eqn:Xwemb} identify the attracting set in question with
\[
\mathring{X}_\Psi
:=
\sX_\Psi(w_0)\cap\prod_{i=1}^n\mathfrak A_i.
\]
Thus $\mathring{X}_\Psi$ is open and dense in the integral variety
$\sX_\Psi(w_0)$.

By the Bia{\l}ynicki-Birula theorem~\cite{BB73}, the attracting map
\[
\mathring{X}_\Psi\longrightarrow X
\]
is an affine-space bundle. It is $G$-equivariant because the actions of
$G$ and $\Gm$ commute. Let $p=B/B\in X$, and let $E$ be the negative
$\Gm$-weight subspace in the decomposition
\begin{equation}\label{eqn:Xtan}
T_p\sX_\Psi(w_0)=T_pX\oplus E.
\end{equation}
Both summands are $B$-stable.

We compute $E$ from the successive projective-space bundles in
Proposition~\ref{prop:idY(k)}. For $e_1(\Psi)\le k\le n$, the relevant
homogeneous space is
\begin{equation}\label{eqn:piece}
G(k)/\bigl(G(k)\cap\tbP(k+1)\bigr)
\cong
\P\bigl(N^\Psi(\varpi_k;k)\bigr)
\subset\P(L(\La_k)).
\end{equation}
At its highest-weight point, the tangent weights are
\begin{equation}\label{eqn:kth}
T_{[\bv_{\La_k}]}
\P\bigl(N^\Psi(\varpi_k;k)\bigr)
\cong
\bigoplus_{k<s\le n}\C_{\epsilon_s-\epsilon_k}
\oplus
\bigoplus_{1\le t\le\hs_k(\Psi)}
\C_{\epsilon_t-\epsilon_k-\delta}.
\end{equation}
Let $\Pi_k$ be the set of weights in~\eqref{eqn:kth}, and put
\begin{equation}\label{eqn:extrawt}
\Pi_k^-:=
\left\{
\epsilon_t-\epsilon_k-\delta
\middle|
1\le t\le\hs_k(\Psi)
\right\}.
\end{equation}
Together with the tangent directions coming from the base $G/Q$, the
classical weights in the first summand of~\eqref{eqn:kth} account for
$T_pX$. Hence the weights of $E$, with multiplicities, form the
disjoint union
\[
\bigsqcup_{k=e_1(\Psi)}^n\Pi_k^-.
\]
By Lemma~\ref{lem:transp}, this is precisely the $\tT$-character of
$\gn(\Psi)\otimes\C_{-\delta}$.

For $e_1(\Psi)\le k \le n$ and $\beta\in\Pi_{k}$, let $U_\beta$ be the one-dimensional root
subgroup of $\widetilde G(\!(z)\!)$ with
$\Lie U_\beta\cong\C_\beta$. 
For $e_1(\Psi) \le k \le n$ and $\beta, \beta' \in \Pi_{k}$, we have
\[
\beta + \beta' \equiv \epsilon_t+\epsilon_{s} - 2 \epsilon_{k} \mod \Z \delta \qquad (t,s \neq k),
\]
which is neither an affine root nor a null-root. It follows that
\begin{equation}
U_{\beta}U_{\beta'} = U _{\beta'} U_{\beta} \subset \widetilde G(\!(z)\!).\label{eqn:Ucomm}
\end{equation}

For each $e_1(\Psi)\le k\le n$, the standard big cells in the
projective homogeneous spaces in~\eqref{eqn:piece}, together with the
successive bundle presentation of Proposition~\ref{prop:idY(k)}, give
an open immersion
\begin{equation}\label{eqn:Umult}
\prod_{k'=k}^n\A^{|\Pi_{k'}|}
\cong
\left(\prod_{\beta\in\Pi_k}U_\beta\right)
\left(\prod_{\beta\in\Pi_{k+1}}U_\beta\right)
\cdots
\left(\prod_{\beta\in\Pi_n}U_\beta\right)\cdot p
\hookrightarrow\sY_\Psi(k),
\end{equation}
where each of $\prod_{\beta\in\Pi_{k'}}U_\beta$ is well-defined by~\eqref{eqn:Ucomm}. Notice that
\[
\sum_{k'=k}^n|\Pi_{k'}|=\dim\sY_\Psi(k).
\]

Together, the big cell in $G/B$ and~\eqref{eqn:Umult}, with
$k=e_1(\Psi)$, give a $\tT$-equivariant surjection
\[
\bigoplus_{i>j} \C E_{ij}
\oplus
\bigl(\gn(\Psi)\otimes\C_{-\delta}\bigr)
\longrightarrow
T_p\sX_\Psi(w_0).
\]
The two sides have the same dimension, so this map is an isomorphism.
In particular,
\[
E\cong\gn(\Psi)\otimes\C_{-\delta}
\]
as $\tT$-modules.

Consider the morphism
\[
\mathsf{exp}_\Psi\colon
\gn(\Psi)\longrightarrow\P(L(\La_n)),
\qquad
\xi \longmapsto
\bigl[\exp( \xi z^{-1})\bv_{\La_n}\bigr].
\]
Since $\gn(\Psi)$ is $B$-stable and
$\C\bv_{\La_n}$ is a $G$-stable line, this morphism is
$B$-equivariant. The degree-$(-1)$ component of the level-one fundamental module
$L(\wp)$ of affine type $\mathsf{A}_{n-1}^{(1)}$ is the adjoint
$G$-module $\mathfrak{sl}(n)$. Hence, under the identification
\[
\mathfrak{sl}(n) \otimes \C_{\varpi_n} \cong L(\La_n)_{-1} \subset L(\La_n) = L (\wp) \otimes \C_{\varpi_n},
\]
the degree-$(-1)$ component of
$\exp(\xi z^{-1})\bv_{\La_n}$ recovers $\xi$. Thus
$\mathsf{exp}_\Psi$ is an isomorphism onto its image.
The root-subgroup chart~\eqref{eqn:Umult} shows that this image is
precisely the attracting fiber over $p$. Therefore the attracting fiber over
$p$ is $B$-equivariantly isomorphic to
$\gn(\Psi)\otimes\C_{-\delta}$.

Since $X\cong G/B$ and the attracting map is $G$-equivariant, we
conclude that
\[
\mathring{X}_\Psi
\cong
G\times^B\bigl(\gn(\Psi)\otimes\C_{-\delta}\bigr)
=T^*_\Psi X.\qedhere
\]
\end{proof}

\begin{cor}[Corollary of the proof of Theorem~\ref{thm:XTid}]
\label{cor:injfib}
The projection of~\eqref{eqn:Xwemb} to $\P(L(\La_n))$ is injective
on every fiber of $T^*_\Psi X\to X$.
\end{cor}

\begin{proof}
This follows from the injectivity of $\mathsf{exp}_\Psi$ at the base
point and from $G$-equivariance.
\end{proof}

By comparison with Lusztig~\cite{Lus81b}, we obtain the following:

\begin{cor}[Ng\^o~\cite{Ngo99}, Mirkovi\'c--Vybornov~\cite{MV03}]\label{cor:Lus}
The composition
\[
\sX_{\Delta^+}(w_0) \hookrightarrow \prod_{i =1}^n \P(L(\La_i)) \rightarrow \P(L(\La_n))
\]
defines a resolution of a compactification of the nilpotent cone of $\mathfrak{gl}(n, \C)$, realized in the affine Grassmannian of $G$.
\end{cor}

For $\la\in\sP$, let $\cO_{T^*_\Psi X}(\la)$ denote the restriction of $\cO_{\sX_\Psi(w_0)}(\la)$ to the open subset in Theorem~\ref{thm:XTid}.

\begin{cor}\label{cor:lines}
Let $\Psi\subset\Delta^+$ be a root ideal, and let $\la\in\sP$. Then
\[
\cO_{T^*_\Psi X}(\la)
\cong
\pi_\Psi^*\cO_X(\la).
\]
\end{cor}

\begin{proof}
Both line bundles are $G$-equivariant, so it suffices to compare their
restrictions to the fiber over $B/B\in X$ as
$(B\times\Gm)$-equivariant line bundles. An equivariant line bundle on
the affine $B$-module $\gn(\Psi)$ is determined by its fiber at the
$\tT$-fixed origin. The two fibers have the same $\tT$-weight by
construction, proving the assertion.
\end{proof}

\begin{cor}\label{cor:vbsub}
Let $\Psi\subset\Delta^+$ be a root ideal, and let $w\in\Sym_n$ be
$\Psi$-tame. Then the $\Gm$-attracting set of
$X(w)\subset\sX_\Psi(w)$ is isomorphic to $T^*_\Psi X(w)$.
\end{cor}

\begin{proof}
By~\eqref{eqn:Xw-associated},
\[
\sX_\Psi(w)
\cong
X_Q(w)\times_{G/Q}\sX_\Psi(w_0).
\]
The attracting set is therefore the restriction of the attracting
bundle in Theorem~\ref{thm:XTid} to
$X(w)=X_Q(w)\times_{G/Q}G/B$, which is precisely $T^*_\Psi X(w)$.
\end{proof}

We also record the nef cone of $\sX_\Psi(w_0)$; see
\cite[Definition~1.4.1]{Laz04}.

\begin{cor}\label{cor:ample}
For any nonempty root ideal $\Psi \subset \Delta^+$, we have $\Pic \sX_\Psi(w_0) \cong \mathsf{P}$. For each $\la \in \mathsf{P}$, the line bundle $\cO_{\sX_\Psi(w_0)}(\la)$ is nef if and only if $\la \in \Par$.
\end{cor}

\begin{proof}
The standard tower of forgetful maps implies that $G/Q$ is a $d_1(\Psi)$-fold iterated projective-space bundle.
Thus Proposition~\ref{prop:idY(k)} and Theorem~\ref{thm:str} present
$\sX_\Psi(w_0)$ as an $n$-fold iterated projective-space bundle. At the
$i$-th stage, $\cO_{\sX_\Psi(w_0)}(\varpi_i)$ restricts to the primitive
ample generator on the corresponding projective-space fiber. Repeated
applications of~\cite[II, Exercise~7.9]{Har77} therefore give $\Pic\sX_\Psi(w_0)\cong\sP$.

If $\cO_{\sX_\Psi(w_0)}(\la)$ is nef, then its restriction to
$X=G/B$ is nef. Hence
\[
\la_i-\la_{i+1}\ge0
\qquad (1\le i<n).
\]
Since $\Psi$ is nonempty, $\hs_n(\Psi)>0$, and the last projective-space
stage contains the line
\[
Y:=
\P\left(
\C\bv_{\La_n}
\oplus
\C(E_{1,n}z^{-1})\bv_{\La_n}
\right)
\cong\P^1
\subset\sY_\Psi(n)\subset\sX_\Psi(w_0).
\]
The construction of the line bundles gives $\cO_{\sX_\Psi(w_0)}(\la)|_Y\cong\cO_{\P^1}(\la_n)$.
Nefness therefore implies $\la_n\ge0$, and consequently
$\la\in\Par$.

Conversely, each $\cO_{\sX_\Psi(w_0)}(\varpi_i)$ is the pullback of
$\cO(1)$ under the $i$-th component of~\eqref{eqn:Xwemb}, and is
therefore nef. Hence $\cO_{\sX_\Psi(w_0)}(\la)$ is nef for every
$\la\in\Par$.
\end{proof}

\section{Consequences}

We retain the notation and assumptions of the preceding sections.

\subsection{Vanishing theorems}

\begin{thm}\label{thm:incl}
Let $\Psi\subset\Delta^+$ be a root ideal. There exists an effective
$(G\times\Gm)$-equivariant Cartier divisor $\partial$ on
$\sX_\Psi(w_0)$ such that
\[
\cO_{\sX_\Psi(w_0)}(\partial)
\cong
\cO_{\sX_\Psi(w_0)}(\varpi_n)\qquad\text{and}\qquad
\operatorname{supp}\partial
=
\sX_\Psi(w_0)\setminus T^*_\Psi X.
\]
Moreover, for every $\Psi$-tame element $w\in\Sym_n$, every
$\la\in\Par$, and every $i\in\Z_{\ge0}$, we have
\begin{equation}\label{eqn:limw}
H^i\bigl(T^*_\Psi X(w),\cO_{T^*_\Psi X(w)}(\la)\bigr)
\cong
\varinjlim_m
H^i\bigl(\sX_\Psi(w),
\cO_{\sX_\Psi(w)}(\la+m\partial)\bigr)
\otimes\C_{m\varpi_n}.
\end{equation}
In particular,
\begin{equation}\label{eqn:CH}
H^{>0}\bigl(T^*_\Psi X(w),
\cO_{T^*_\Psi X(w)}(\la)\bigr)=0.
\end{equation}
\end{thm}

\begin{proof}
Consider the last component of the closed immersion~\eqref{eqn:Xwemb}:
\[
\sX_\Psi(w_0)
\longrightarrow
\P\bigl(H^0\bigl(\sX_\Psi(w_0),\cO(\varpi_n)\bigr)^\vee\bigr)
\subset
\P(L(\La_n)).
\]
Let $\tau$ be the section of $\cO(\varpi_n)$ dual to the
highest-weight vector $\bv_{\La_n}$, and let $\partial$ be its zero
divisor. Since $\sX_\Psi(w_0)$ is integral and $\tau$ is nonzero,
$\partial$ is an effective Cartier divisor.

The description of the attracting fiber over $B/B$ in the proof of
Theorem~\ref{thm:XTid} identifies its last component with
\[
\xi \longmapsto
\bigl[\exp(\xi z^{-1})\bv_{\La_n}\bigr]
\qquad (\xi \in\gn(\Psi)).
\]
The coefficient of $\bv_{\La_n}$ in this expression is~$1$.
Consequently, $\tau$ does not vanish on the attracting fiber over
$B/B$, and $G$-equivariance gives
\[
\operatorname{supp}\partial\cap T^*_\Psi X=\varnothing.
\]
Conversely, in the affine chart $\{\tau\ne0\}$, the degree-$(-1)$
component of $\exp(\xi z^{-1})\bv_{\La_n}$ recovers $\xi$. Thus the image
of the fiber over $B/B$ is closed in this affine chart, and its
boundary lies in the hyperplane $\{\tau=0\}$. Taking the
$G$-translates of this fiber, we obtain
\[
\operatorname{supp}\partial
=
\sX_\Psi(w_0)\setminus T^*_\Psi X.
\]

For a $\Psi$-tame element $w$, Corollary~\ref{cor:vbsub} identifies
$T^*_\Psi X(w)$ with the intersection of $\sX_\Psi(w)$ and the
principal open subset $\{\tau\ne0\}$. Thus the open immersion
\[
\jmath\colon T^*_\Psi X(w)\hookrightarrow\sX_\Psi(w)
\]
is affine. Corollary~\ref{cor:lines} gives
\[
\cO_{T^*_\Psi X(w)}(\la)
\cong
\pi_\Psi^*\cO_{X(w)}(\la),
\]
and localization with respect to $\tau$ yields
\begin{equation}\label{eqn:jstar-localization}
\jmath_*\cO_{T^*_\Psi X(w)}(\la)
\cong
\varinjlim_m
\cO_{\sX_\Psi(w)}(\la+m\partial)
\otimes\C_{m\varpi_n}.
\end{equation}
Here the transition map from the $m$-th term to the $(m+1)$-st is
multiplication by $\tau$; the character twist makes this map
$\tT$-equivariant.

Since $\jmath$ is affine, cohomology on $T^*_\Psi X(w)$ is the
cohomology of the sheaf on the left-hand side of
\eqref{eqn:jstar-localization}. Compatibility of cohomology with
direct limits~\cite[III, Proposition~2.9]{Har77} therefore gives
\eqref{eqn:limw}. Finally, since $\la+m\varpi_n\in\Par$ for every $m\ge0$,
Theorem~\ref{thm:Xmain}(1) gives~\eqref{eqn:CH}.
\end{proof}

\begin{rem}\label{rem:coverage}
The vanishing result~\eqref{eqn:CH} establishes the tame case of the vanishing conjecture proposed by Blasiak--Morse--Pun~\cite[Conjecture 3.4(ii)]{BMP}. Since $w_0$ is
$\Psi$-tame for every root ideal $\Psi$, the theorem applies to all
ordinary Catalan functions $H(\Psi;\la)$, including their
$k$-Schur specializations; see Remark~\ref{rem:rotation-conventions}. In particular,~\eqref{eqn:CH} proves the vanishing conjectures of Chen--Haiman~\cite[Conjecture 5.4.3(2)]{Che10} and Shimozono--Weyman~\cite[Conjecture 5]{SW00}. This result was previously known in the case where $\la$ is strictly dominant~\cite{Pan10,MvdK92}, or in certain special cases~\cite{Bro93,Bro94,Hag09}. However, these results do not cover all cases in which $H(\Psi; \la; w_0)$ is a $k$-Schur polynomial~\cite{BMPS}, or all cases in which $\gn(\Psi)$ is the nilpotent radical of a proper parabolic subgroup of~$G$.
\end{rem}

\begin{cor}[{\cite[Conjecture 3.4(iv)]{BMP}}]\label{cor:BMP}
Let $\Psi \subset \Delta^+$ be a root ideal, let $w \in \Sym_n$ be $\Psi$-tame, and let $\la \in \Par$. Then
\[
H^{0} \bigl( T^*_\Psi X (w), \cO_{T^*_\Psi X (w)} ( \la ) \bigr)
\]
admits an excellent filtration in the sense of van der Kallen~\cite{vdK89}.
\end{cor}

\begin{proof}
By Theorem~\ref{thm:Xmain}(3), every $B$-module in the direct system on the right-hand side of~\eqref{eqn:limw} admits an excellent filtration. Moreover, by~\cite[Corollary~1.8]{vdK89}, the inductive limit of such modules also admits an excellent filtration. The claim follows from~\eqref{eqn:limw} with $i=0$.
\end{proof}

\begin{cor}[Parabolic vanishing]\label{cor:pv}
Let $P \subset G$ be a parabolic subgroup containing $B$, and let $\Psi \subset \Delta^+$ be a root ideal such that $\mathfrak{n}(\Psi)$ is $P$-stable. Define
\[
T^*_\Psi X^P := G \times^P \mathfrak{n}(\Psi) \stackrel{\pi^P_\Psi}{\longrightarrow} G/P =: X^P.
\]
Suppose that $\la\in\Par$ satisfies
\[
\langle\alpha_i,\la\rangle=0
\qquad\text{whenever }P_i\subset P.
\]
Then the line bundle $\cO_X(\la)$ on $X = G/B$ descends to a line bundle $\cO_{X^P}(\la)$ on $X^P$, and we have
\[
H^{>0}\bigl(T^*_\Psi X^P, (\pi^P_\Psi)^* \cO_{X^P}(\la)\bigr) = 0.
\]
\end{cor}

\begin{proof}
Since $P$ stabilizes $\mathfrak{n}(\Psi)$, there is a natural $P/B$-fibration
\[
\eta : T^*_\Psi X = G \times^B \mathfrak{n}(\Psi) \longrightarrow G \times^P \mathfrak{n}(\Psi) = T^*_\Psi X^P.
\]
The condition on $\la$ ensures that the weight $\la$ descends to a character of $P$, and hence the line bundle $\cO_X(\la)$ is the pullback of the $G$-equivariant line bundle $\cO_{X^P}(\la)$ on $G/P$. In particular, $(\pi_\Psi)^*\cO_X(\la)$ is trivial along the fibers of $\eta$.

By the Borel--Weil--Bott theorem, we have $H^p(P/B, \cO_{P/B}) = \C^{\delta_{p,0}}$. By flat base change~\cite[I\!I\!I, Proposition 9.3]{Har77}, we deduce $\mathbb{R}^p \eta_* (\pi_\Psi)^* \cO_X(\la) = 0$ for all $p > 0$. Hence, the Leray spectral sequence
\[
E_2^{q,p} := H^{q} \bigl( T^*_\Psi X^P, \mathbb{R}^{p} \eta_* ( \pi_\Psi )^* \cO_X(\la) \bigr) \Rightarrow H^{q+p} \bigl( T^*_\Psi X, ( \pi_\Psi )^* \cO_X(\la) \bigr)
\]
degenerates at the $E_2$-page. This yields
\[
E_2^{q,0} \cong H^{q} \bigl( T^*_\Psi X^P, \eta_* ( \pi_\Psi )^* \cO_X(\la) \bigr) = H^{q} \bigl( T^*_\Psi X^P, (\pi^P_\Psi)^* \cO_{X^P}(\la) \bigr).
\]
Therefore,~\eqref{eqn:CH} for $q > 0$ and $w = w_0$ implies the claimed result.
\end{proof}

\begin{rem}\label{rem:gen}
\textbf{(1)} Corollary~\ref{cor:pv} has a $B$-equivariant analogue, in the same sense as~\eqref{eqn:CH}. \textbf{(2)} The results of~\S\ref{subsec:fv} are valid over an arbitrary base field. The results of~\S\ref{subsec:adm} also remain valid in positive characteristic~\cite{Jos06}, except when $n = 2$, where the corresponding affine Lie algebra is not simply-laced. Consequently, all the results in~\S\ref{sec:XP}, as well as Theorem~\ref{thm:incl} and Corollary~\ref{cor:pv}, remain valid in arbitrary characteristic for $n \ge 3$. The case $n = 2$ in positive characteristic can be handled separately by elementary arguments.
\end{rem}

\subsection{Simple head property}\label{subsec:sh}

\begin{lem}\label{lem:infiinj}
Let $\Psi\subset\Delta^+$ be a root ideal. Then $T^*_\Psi X$ admits
a natural infinitesimal action of $\mathfrak{gl}(n,\C[z])$. For every
$\la\in\Par$, this action equips
\[
H^0\bigl(T^*_\Psi X,\cO_{T^*_\Psi X}(\la)\bigr)
\]
with the structure of a graded $\mathfrak{gl}(n,\C[z])$-module, and
the restriction map
\[
H^0\bigl(\sX_\Psi(w_0),\cO_{\sX_\Psi(w_0)}(\la)\bigr)
\longrightarrow
H^0\bigl(T^*_\Psi X,\cO_{T^*_\Psi X}(\la)\bigr)
\]
is an injective homomorphism of graded
$\mathfrak{gl}(n,\C[z])$-modules.
\end{lem}

\begin{proof}
Differentiating the $\tbG$-action and restricting it to
$\C\mathrm{Id}\oplus\mathfrak{sl}(n,\C[z])\subset\Lie\tbG$, we obtain
a $\mathfrak{gl}(n,\C[z])$-action by letting
$z\C[z]\mathrm{Id}$ act trivially. The loop-rotation action makes
this action graded.

The restriction map is injective because $\sX_\Psi(w_0)$ is integral
and $T^*_\Psi X$ is dense in it. Since the action on $T^*_\Psi X$ is
obtained by restricting that on $\sX_\Psi(w_0)$ to the open subset,
the restriction map is compatible with the Lie algebra actions.
\end{proof}

\begin{rem}
We caution that the $\gl(n, \C[z])$-action on
\[
H^{0} \bigl( T^*_{\Psi} X, \cO_{T^*_{\Psi} X} ( \la ) \bigr) \qquad (\la \in \Par)
\]
is, in general, \emph{not} compatible with the identification~\eqref{eqn:limw}.
This is
analogous to the fact that the $\mathfrak g$-module inclusion
\[
H^0 \bigl( X , \cO_{X} ( \la ) \bigr) \hookrightarrow H^0 \bigl( w_0 B w_0 B/B, \cO_{X} ( \la ) \bigr) \qquad (\la \in \Par)
\]
is not compatible with character twists when regarded as a map of $\mathfrak{b}$-modules.
\end{rem}

\begin{thm}\label{thm:infi}
For each $\la \in \Par$ and each root ideal $\Psi \subset \Delta^+$, the $\tbG$-module
\[
H^{0} \bigl( \sX_{\Psi} ( w_0 ), \cO_{\sX_\Psi ( w_0 )} ( \la ) \bigr)
\]
has a simple head isomorphic to $H^0 \bigl( X, \cO_X ( \la ) \bigr) \cong V (\la)^{\vee}$.
\end{thm}

We first record a consequence of Theorem~\ref{thm:infi}.

\begin{cor}\label{cor:infi}
Let $\Psi \subset \Delta^+$ be a root ideal, let $w \in \Sym_n$ be $\Psi$-tame, and let $\la \in \Par$. Then the $\tbB$-module
\[
H^{0} \bigl( \sX_{\Psi} ( w ), \cO_{\sX_\Psi ( w )} ( \la ) \bigr)
\]
has a simple head.
\end{cor}

\begin{proof}
Choose a reduced expression for $w$ and extend it, by prepending
finite simple reflections, to a reduced expression for $w_0$. By repeated applications of Corollary~\ref{cor:D-branch} to the presentations in~\eqref{eqn:defN}, we obtain a surjective $\tbB$-module map
\begin{equation}
H^{0} \bigl( \sX_{\Psi} ( w_0 ), \cO_{\sX_\Psi ( w_0 )} ( \la ) \bigr) \longrightarrow H^{0} \bigl( \sX_{\Psi} ( w ), \cO_{\sX_\Psi ( w )} ( \la ) \bigr). \label{eqn:Hcov}
\end{equation}
By Theorem~\ref{thm:infi} and the PBW theorem, the source
of~\eqref{eqn:Hcov} has a simple head as a $\tbB$-module.
Since the target is a nonzero quotient of this module, it also has
a simple head.
\end{proof}

We now prepare the proof of Theorem~\ref{thm:infi}. 
For each $2 \le r \le n$, we define Lie subalgebras
\begin{align*}
\gp_{(r)} &:= \Span \{ E_{1,1}, E_{i,j} \mid (1 \le i \le r,\ 2 \le j \le r)\} \subset \gl(n), \\
\gp_{(r)}^- &:= \Span \{ E_{1,1}, E_{i,j} \mid (2 \le i \le r,\ 1 \le j \le r) \} \subset \gl(n).
\end{align*}
We then set
\[
\g_{(r)} := \gp_{(r)} + \gp_{(r)}^-, \qquad \gs_{(r)} := \gp_{(r)} \cap \gp_{(r)}^-,
\]
so that
\[
\gl(r) = \g_{(r)} \subset \gl(n) \supset \gs_{(r)} = \C \oplus \gl(r-1) .
\]
We also introduce
\[
\tp_{(r)} := \gp_{(r)} + \tb
\subset \tg_{(r)} := \g_{(r)} + \tb
\subset \Lie \tbG.
\]

\begin{prop}\label{prop:key}
Fix $2 \le r \le n$, and let $M$ be a finite-dimensional $\tp_{(r)}$-module which is semisimple as a $\tT$-module (with respect to the integrated $\wgt$-action). Then there exists a surjection of $\gp_{(r)}$-modules
\[
U ( \gp_{(r)}) \otimes_{U ( \gs_{(r)})} M \twoheadrightarrow 
\mathscr D_{s_{r-1} s_{r-2} \cdots s_1} ( M^{\vee} )^{\vee}.
\]
\end{prop}

\begin{proof}
Let $\tbP_{(r)}$ and $\tbG_{(r)}$ be the proalgebraic subgroups of
$\tbG$ with Lie algebras $\tp_{(r)}$ and $\tg_{(r)}$, respectively,
and let $P_{(r)}\subset G_{(r)}$ be the connected algebraic subgroups
of $G$ with Lie algebras $\gp_{(r)}\subset\g_{(r)}$. Multiplication
induces
\[
G_{(r)}/P_{(r)}
\cong
\tbG_{(r)}/\tbP_{(r)}
\cong
\P^{r-1}.
\]

Since $M$ is finite-dimensional, it acquires the structure of a rational $\tbP_{(r)}$-module. 
After a cyclic relabeling of the affine Dynkin diagram,
Lemma~\ref{lem:ind} gives
\begin{equation}
\mathscr D_{s_{r-1} s_{r-2} \cdots s_1} ( M^{\vee} )^{\vee} 
\cong H^0 \bigl( \P^{r-1}, \cE ( M )\bigr),\label{eqn:infl-small}
\end{equation}
where $\cE(M)$ denotes the sheaf of sections of $\tbG_{(r)}\times^{\tbP_{(r)}}M^\vee$.

Let $U_{(r)}^-$ be the opposite unipotent radical of
$P_{(r)}\subset G_{(r)}$, with Lie algebra $\gu_{(r)}^-$. Restriction
to the big cell and evaluation at its identity give an
$\gs_{(r)}$-equivariant map
\begin{equation}\label{eqn:comploc}
\imath\colon
H^0\bigl(\P^{r-1},\cE(M)\bigr)
\hookrightarrow
\C[U_{(r)}^-]\otimes M
\longrightarrow M.
\end{equation}
Its restriction to
$H^0\bigl(\P^{r-1},\cE(M)\bigr)^{\gu_{(r)}^-}$ is injective: a $U_{(r)}^-$-invariant
section is constant on the big cell and is therefore determined by
its value at the identity.

Since $H^0(\P^{r-1},\cE(M))$ is a finite-dimensional
$\g_{(r)}$-module, the $\gs_{(r)}$-module
$H^0\bigl(\P^{r-1},\cE(M)\bigr)^{\gu_{(r)}^-}$ is semisimple and generates $H^0(\P^{r-1},\cE(M))$ under the
$\gp_{(r)}$-action. The injection
$H^0(\P^{r-1},\cE(M))^{\gu_{(r)}^-}\hookrightarrow M$ in~\eqref{eqn:comploc} splits as
a homomorphism of $\gs_{(r)}$-modules because $\gs_{(r)}$ is
reductive. Induction from $\gs_{(r)}$ to $\gp_{(r)}$ therefore gives
a surjection
\[
U(\gp_{(r)})\otimes_{U(\gs_{(r)})}M
\twoheadrightarrow
U(\gp_{(r)})\otimes_{U(\gs_{(r)})}H^0(\P^{r-1},\cE(M))^{\gu_{(r)}^-}
\twoheadrightarrow H^0\bigl(\P^{r-1},\cE(M)\bigr).
\]
Together with~\eqref{eqn:infl-small}, this proves the proposition.
\end{proof}

\begin{proof}[Proof of Theorem~\ref{thm:infi}]
For $e_1(\Psi) \le k \le n$, let $\gu(k)$ denote the Lie algebra of the unipotent radical of the subgroup $G(k)\cap G(k+1)\tbB$ inside $G(k)$; see Lemma~\ref{lem:Ustr}.

Let $\Phi$ be the affine
Dynkin diagram automorphism that increases every index in $\tI_\af$ by one
modulo~$n$. Then
\[
\Phi^{k-1}\bigl(\gp_{(n-k+\hs_k(\Psi)+1)}\bigr) = \g(k)\cap \bigl( \g(k+1) + \Lie \tbB \bigr)
\qquad (e_1(\Psi) \le k \le n).
\]
For the remaining finite stages, put
\[
\gu(k) := \bigoplus_{j=k+1}^{n} \C E_{kj} \subset \g \qquad (1 \le k \le d_1(\Psi)).
\]
Then $\gu(1)\oplus\cdots\oplus\gu(d_1(\Psi))$ is the nilpotent
radical of the parabolic subalgebra
$\g\cap\Lie\tbP(e_1(\Psi))$. Its $\tT$-weights are
\begin{equation}\label{eqn:T0wt}
\epsilon_i-\epsilon_j
\qquad(1\le i\le d_1(\Psi),\ i<j\le n).
\end{equation}

We expand as
\[
\la = \sum_{i=1}^n m_i ( \la ) \varpi_i \qquad\text{and set} \qquad\La_{(k)} := \sum_{j=1}^k m_j ( \la ) \La_j \qquad (0 \le k \le n).
\]

Applying Proposition~\ref{prop:key}, after the
cyclic relabeling $\Phi^{1-k}$, to the successive associated-bundle construction of $\sY_\Psi(k)$ in Proposition~\ref{prop:idY(k)} gives a surjection
\begin{align}\nonumber
U(\gu(k))\otimes
H^0\bigl(\sY_\Psi(k+1),\cO_{\sY_\Psi(k+1)}(\la)\bigr)
\otimes\C_{-\La_{(k)}}\!&\,\\
\twoheadrightarrow
H^0\bigl(\sY_\Psi(k),\cO_{\sY_\Psi(k)}(\la)\bigr)\!&\,
\otimes\C_{-\La_{(k-1)}}
\label{eqn:Y-cyclic-surj}
\end{align}
of $U(\gu(k))$-modules for $e_1(\Psi)\le k\le n$, semisimple with respect to the $\tT$-action.

Iterating~\eqref{eqn:Y-cyclic-surj}
therefore gives a surjection
\begin{align}\nonumber
U(\gu(e_1(\Psi)))\otimes\cdots\otimes U(\gu(n))\!&\,
\otimes\C_{-\La_{(n)}}\\
\twoheadrightarrow \!&\,\,
H^0\bigl(\sX_\Psi(w_0^\Psi),
\cO_{\sX_\Psi(w_0^\Psi)}(\la)\bigr)
\otimes\C_{-\La_{(d_1(\Psi))}}.\label{eqn:Y-PBW-surj}
\end{align}

For $e_1(\Psi)\le k\le n$, the $\tT$-weights of $\gu(k)$ are
\begin{equation}\label{eqn:wtuk}
\epsilon_k - \epsilon_j \quad (k < j \le n), 
\qquad 
\epsilon_k - \epsilon_i + \delta \quad (1 \le i \le \hs_k(\Psi))
\end{equation}
(cf.~\eqref{eqn:extrawt}). These weights are pairwise distinct as $k$ varies. They are also
disjoint from the weights in~\eqref{eqn:T0wt}.

The Demazure functors needed to pass from $w_0^\Psi$ to $w_0$
can be written as
\[
\mathscr D_{w_0w_0^\Psi}
\cong
\mathscr D_{s_{n-1}\cdots s_1}
\circ\mathscr D_{s_{n-1}\cdots s_2}
\circ\cdots\circ
\mathscr D_{s_{n-1}\cdots s_{d_1(\Psi)}}.
\]
Applying Proposition~\ref{prop:key} to these 
stages and combining the result with~\eqref{eqn:Y-PBW-surj} gives
\begin{equation}\label{eqn:X-PBW-surj}
U(\gu(1))\otimes U(\gu(2))\otimes\cdots\otimes U(\gu(n))
\otimes\C_{-\La_{(n)}}
\twoheadrightarrow
H^0\bigl(\sX_\Psi(w_0),
\cO_{\sX_\Psi(w_0)}(\la)\bigr).
\end{equation}

The weight descriptions~\eqref{eqn:T0wt} and~\eqref{eqn:wtuk}
show that
\[
\gu(1)+\cdots+\gu(n)
=\bigoplus_{i=1}^n\gu(i)\subset\Lie\tbG
\]
is multiplicity-free as a $\tT$-module. The PBW theorem and~\eqref{eqn:X-PBW-surj} therefore imply that $H^0\bigl( \sX_\Psi(w_0), \cO_{\sX_\Psi(w_0)}(\la) \bigr)$ is generated by
the image of the one-dimensional tensor factor
$\C_{-\La_{(n)}}$ as a $\Lie\tbG$-module. This cyclic line has rotation degree zero and lies
in
\[
H^0\bigl(X,\cO_X(\la)\bigr)\cong V(\la)^{\vee},
\]
which completes the proof.
\end{proof}

\subsection{Monotonicity of multiplicities}\label{subsec:mon}

\begin{prop}\label{prop:Psi-surj}
Let $\Psi' \subset \Psi \subset \Delta^+$ be root ideals, and let $w', w \in \Sym_n$ be $\Psi$-tame elements such that $X(w') \subset X(w)$. Then, for every $\la \in \Par$, there is an inclusion of $\tbB$-modules
\[
N_{w'}^{\Psi'}(\la) \subset N_w^{\Psi}(\la).
\]
\end{prop}

\begin{proof}
Since $d_1(\Psi')\ge d_1(\Psi)$, every $\Psi$-tame element is
$\Psi'$-tame. In particular, $w'$ is $\Psi'$-tame.

Choose the sequence $\bi=\bi(\Psi,w)$ from~\eqref{eqn:defbiPw} so
that its finite part contains, by the subword property for the
Bruhat order, a reduced expression for $w'$. The inclusion $\Psi'\subset\Psi$ implies
$\hs_k(\Psi')\le\hs_k(\Psi)$ for each
$1\le k\le n$. Hence the chosen reduced expression for
$w'$, followed by the blocks defining $N_{w'}^{\Psi'}(\la)$, forms
a subsequence $\bi'\subset\bi$ that realizes
$N_{w'}^{\Psi'}(\la)$. Lemma~\ref{lem:demseq}(2) gives a closed
immersion
\[
Z(\bi')\hookrightarrow Z(\bi),
\]
and the line bundle $\mathcal L_w^\Psi(\la)$ restricts to
$\mathcal L_{w'}^{\Psi'}(\la)$. Thus~\eqref{eqn:N-BS-sections}
gives a restriction homomorphism
\begin{equation}\label{eqn:Nrest}
N_w^\Psi(\la)^\vee
\cong
H^0\bigl(Z(\bi),\mathcal L_w^\Psi(\la)\bigr)
\longrightarrow
H^0\bigl(Z(\bi'),\mathcal L_{w'}^{\Psi'}(\la)\bigr)
\cong
N_{w'}^{\Psi'}(\la)^\vee.
\end{equation}

The $\tT$-weights of the simple heads on both sides
of~\eqref{eqn:Nrest}, as described in Corollary~\ref{cor:infi},
coincide: both are realized as the (dual of the) fiber of
$\mathcal L_w^\Psi(\la)$ at the $\tT$-fixed point
$Z(\varnothing)\subset Z(\bi')\subset Z(\bi)$. Under this
identification, the image of~\eqref{eqn:Nrest} maps nontrivially to the
simple head of $N_{w'}^{\Psi'}(\la)^\vee$. Hence~\eqref{eqn:Nrest} is
surjective. Taking duals yields the desired inclusion of
$\tbB$-modules.
\end{proof}

\begin{prop}\label{prop:incl}
Let $\Psi' \subset \Psi \subset \Delta^+$ be root ideals, and let $w',w \in \Sym_n$ be $\Psi$-tame elements such that $X( w' ) \subset X ( w )$. Then there is a closed immersion $\sX_{\Psi'} (w') \subset \sX_\Psi(w)$ that induces, for every $\la \in \Par$, a surjection
\[
H^0 \bigl( \sX_\Psi(w), \cO_{\sX_\Psi(w)} ( \la )\bigr) \twoheadrightarrow H^0 \bigl( \sX_{\Psi'}(w'), \cO_{\sX_{\Psi'}(w')} ( \la )\bigr).
\] 
\end{prop}

\begin{proof}
Recall that the homogeneous coordinate ring of $\sX_\Psi(w)$ is $\bigoplus_{\la \in \Par} ( N_w^\Psi ( \la ) )^{\vee}$. By Proposition~\ref{prop:Psi-surj}, the natural map
\[
N_w^\Psi ( \la )^{\vee} \longrightarrow N_{w'}^{\Psi'} ( \la )^{\vee}
\]
is surjective for each $\la \in \Par$. These surjections are compatible with multiplication. It follows that the homogeneous coordinate ring of $\sX_{\Psi'}(w')$ is a quotient of that of $\sX_\Psi(w)$. Applying Theorem~\ref{thm:Xmain}(2), we obtain the desired surjection of global sections.
\end{proof}

\begin{cor}\label{cor:SW}
Let $\Psi' \subset \Psi \subset \Delta^+$ be root ideals, and let $w', w \in \Sym_n$ be $\Psi$-tame elements such that $X(w') \subset X(w)$. Then, for all $\la \in \Par$, the natural restriction map
\[
H^0 \bigl( T^*_\Psi X(w), \cO_{T^*_\Psi X(w)} ( \la ) \bigr) \longrightarrow H^0 \bigl( T^*_{\Psi'} X(w'), \cO_{T^*_{\Psi'} X(w')} ( \la ) \bigr)
\]
is surjective. Moreover, there is a scheme-theoretic identification
\[
\sX_{\Psi'}(w') = \overline{T_{\Psi'}^* X(w')} \subset \sX_\Psi(w).
\]
\end{cor}

\begin{proof}
Note that $w'$ is $\Psi'$-tame. For every $m \ge 0$, Theorem~\ref{thm:incl} gives the following commutative diagram:
\[
\xymatrix{
H^0 \bigl( T^*_\Psi X(w), \cO ( \la ) \bigr) \ar[r] &
H^0 \bigl( T^*_{\Psi'} X(w'), \cO ( \la ) \bigr) \\
H^0 \bigl( \sX_\Psi(w), \cO ( \la + m \varpi_n ) \bigr) \otimes \C_{m \varpi_n} \ar@{^{(}->}[u] \ar@{->>}[r] &
H^0 \bigl( \sX_{\Psi'}(w'), \cO ( \la + m \varpi_n ) \bigr) \otimes \C_{m \varpi_n} \ar@{^{(}->}[u]
},
\]
where the surjection in the bottom line is by Proposition~\ref{prop:incl}.
Taking direct limits over $m$ in the above diagram and applying
\eqref{eqn:limw} with $i=0$ proves the first assertion. The second assertion follows from a direct comparison of the homogeneous coordinate rings via the above commutative diagram by varying $\la \in \Par$.
\end{proof}

\begin{defn}
For a root ideal $\Psi \subset \Delta^+$ and dominant weights $\la, \mu \in \sP^+$, we define the graded multiplicity series by
\[
K_{\la,\mu}^{\Psi} (q) := \sum_{m \in \Z} q^m \dim \Hom_{G \times \Gm^\ro} \bigl( V(\la) \boxtimes \C_{-m\delta}, \, H^0 ( T^*_\Psi X, \cO_{T^*_\Psi X} ( \mu ) )^{\vee} \bigr) \in \Z[\![q]\!].
\]
\end{defn}

The following proves and generalizes~\cite[Conjecture~12]{SW00}.

\begin{cor}\label{cor:SWn}
Let $\Psi' \subset \Psi \subset \Delta^+$ be root ideals. Then, for all $\la, \mu \in \sP^+$, we have a coefficientwise inequality
\[
K_{\la,\mu}^{\Psi'} (q) \le K_{\la,\mu}^{\Psi} (q).
\]
\end{cor}

\begin{proof}
Since rational representations of $(G \times \Gm^\ro)$ are completely reducible, $K_{\la,\mu}^{\Psi}(q)$ records the graded multiplicities of $V(\la)$ in
\[
H^0 \bigl( T^*_\Psi X, \cO_{T^*_\Psi X} ( \mu )\bigr)^{\vee}.
\]
Thus, the case $w = w' = w_0$ of Corollary~\ref{cor:SW} implies the desired inequality.
\end{proof}

\begin{rem}\label{rem:SW}
Corollary~\ref{cor:SW} gives a closed immersion
\[
\Spec \, H^0 \bigl( T^*_\Psi X, \cO_{T^*_\Psi X} \bigr) \longrightarrow \Spec \, H^0 \bigl( T^* X, \cO_{T^* X}\bigr) \subset \mathfrak{sl}(n).
\]
Its image is irreducible and reduced\footnote{This property does not hold if $G$ is replaced by a group of a different type, even when considering an equivariant vector subbundle of $T^*(G/B)$ arising from the pullback of $T^*(G/P)$ for a parabolic subgroup $P \subset G$; see, e.g.,~\cite{CM93,Nam06}.}. Hence it is the closure of a
nilpotent orbit, which we denote by $\overline{\mathbb O_\Psi}$.

Since $\varpi_n$ is the determinant character of $G$, we have
\[
\cO_{\sX_\Psi(w_0)}(\varpi_n)\big|_{T^*_\Psi X} \cong \cO_{T^*_\Psi X} \otimes \C_{-\varpi_n}.
\]
It follows that
\[
K_{\la, k\varpi_n}^{\Psi'}(q) \le K_{\la, k\varpi_n}^{\Psi}(q) \hskip 5mm (k \in \Z, \la \in \sP^+),
\]
whenever $\Psi, \Psi' \subset \Delta^+$ satisfy $\mathbb O_{\Psi'} \subset \overline{\mathbb O_\Psi}$; equality holds if $\mathbb O_{\Psi'} = \mathbb O_\Psi$.

In the case where $\gn(\Psi)$ is the Lie algebra of the unipotent radical of a parabolic subgroup of $G$ corresponding to a composition $\nu$, the nilpotent orbit $\mathbb O_\Psi$ is the orbit associated to the transpose of the partition obtained by rearranging $\nu$~\cite[Theorems~7.1.3 and~7.2.3]{CM93}.

This observation appears implicitly in~\cite[\S5.1]{FS21}.
\end{rem}

The following result generalizes Remark~\ref{rem:SW} and yields a
root-ideal extension of~\cite[Conjecture~13]{SW00}.

\begin{cor}\label{cor:SW13}
Let $1 \le a < b \le n$ and let $\mu \in \Par$ satisfy
\[
\mu_a = \mu_{a+1} = \cdots = \mu_b.
\]
Let $\Psi', \Psi \subset \Delta^+$ be two root ideals satisfying:
\begin{enumerate}
\item $E_{a-1,j}\in\gn(\Psi')\cap\gn(\Psi)$ for $j\ge a$ when
$a>1$;
\item $E_{i,b+1}\in\gn(\Psi')\cap\gn(\Psi)$ for $i\le b$
when $b<n$.
\end{enumerate}
Let $G_{a,b} := \SL(b-a+1) \subset G$ denote the $T$-stable subgroup whose set of roots is $\{\epsilon_i-\epsilon_j\mid a\le i,j\le b\}$. If
\begin{equation}
G_{a,b}\cdot \gn ( \Psi' ) \subset \overline{G_{a,b} \cdot \gn ( \Psi )},\label{eqn:action}
\end{equation}
then we have a coefficientwise inequality
\begin{equation}
K_{\la,\mu}^{\Psi'} (q) \le K_{\la,\mu}^{\Psi} (q)\qquad(\la\in\sP^+).\label{eqn:SW13}
\end{equation}
\end{cor}

\begin{proof}
Let $P:=G_{a,b}B\subset G$. The first two assumptions imply that
$\mathfrak g_{a,b}+\gn(\Psi)$ and
$\mathfrak g_{a,b}+\gn(\Psi')$ are $P$-stable. Since the former is
closed and contains $\gn(\Psi)$,~\eqref{eqn:action} yields
\[
\mathfrak g_{a,b}+\gn(\Psi')
\subset
\mathfrak g_{a,b}+\gn(\Psi).
\]
Hence, for $\Xi=\Psi,\Psi'$, we obtain morphisms
\[
f_\Xi\colon T^*_\Xi X=G\times^B\gn(\Xi)
\longrightarrow
G\times^P\bigl(\mathfrak g_{a,b}+\gn(\Psi)\bigr).
\]

Since $\mu$ is constant on the indices $a,\ldots,b$, its character
extends from $B$ to $P$. By the above decomposition and Remark~\ref{rem:SW}, applied to the
root ideals corresponding to
$\gn(\Xi)\cap\mathfrak g_{a,b}$, the inclusion
\[
\overline{G_{a,b}\cdot\gn(\Psi')}
\subset
\overline{G_{a,b}\cdot\gn(\Psi)}
\]
induces a surjection
\[
(f_\Psi)_*\cO_{T^*_\Psi X}(\mu)
\twoheadrightarrow
(f_{\Psi'})_*\cO_{T^*_{\Psi'}X}(\mu).
\]
Taking global sections gives
\begin{equation}
\imath : H^0 \bigl( T^*_{\Psi} X, \cO_{T^*_{\Psi} X} ( \mu ) \bigr) \longrightarrow H^0 \bigl( T^*_{\Psi'} X, \cO_{T^*_{\Psi'} X} ( \mu ) \bigr).\label{eqn:isurj}
\end{equation}

This map fits into a commutative diagram obtained by restricting sheaves on $G \times^{P} ( \mathfrak g_{a,b} + \gn ( \Psi) )$:
\[
\xymatrix{
H^0 \bigl( T^* X, \cO_{T^* X} ( \mu ) \bigr)\ar@{->>}[d] \ar@{->>}[rd]\\
H^0 \bigl( T^*_{\Psi} X, \cO_{T^*_{\Psi} X} ( \mu ) \bigr) \ar[r]^{\imath\quad}& H^0 ( T^*_{\Psi'} X, \cO_{T^*_{\Psi'} X} ( \mu ) )
}.
\]
The two maps from the top term are surjective by
Corollary~\ref{cor:SW}; hence $\imath$ is surjective.
Taking graded characters of~\eqref{eqn:isurj} yields~\eqref{eqn:SW13}, as desired.
\end{proof}

\medskip

{\small
\textbf{Acknowledgment.} This work was partially supported by JSPS KAKENHI Grant Numbers JP19H01782 and JP24K21192.
The author is deeply grateful to Shrawan Kumar for insightful discussions concerning the scope of this work, to Ben Elias for comments that prompted the inclusion of the references~\cite{MvdK92,Ngo99,MV03,Hag09}, and to Mark Shimozono for helpful correspondence regarding Corollary~\ref{cor:SW13}.}

{\footnotesize
\bibliography{ref}
\bibliographystyle{hplain}}
\end{document}